\documentclass[a4paper, 11pt]{article}
\usepackage{amsmath}
\usepackage{amsfonts}
\usepackage{amssymb}
\usepackage[english]{babel}
\usepackage{graphicx}
\usepackage{amsthm}
\usepackage[title]{appendix}

\usepackage{longtable}
\usepackage{tabularx}
\allowdisplaybreaks

\newcommand{\C}{{\cal C}}

\newcommand{\Hh}{{\cal H}}

\newtheorem{theorem}{Theorem}[section]

\newtheorem{lemma}[theorem]{Lemma}

\def\whitebox{{\hbox{\hskip 1pt
 \vrule height 6pt depth 1.5pt
 \lower 1.5pt\vbox to 7.5pt{\hrule width
    3.2pt\vfill\hrule width 3.2pt}%
 \vrule height 6pt depth 1.5pt
 \hskip 1pt } }}
\def\qed{\ifhmode\allowbreak\else\nobreak\fi\hfill\quad\nobreak
     \whitebox\medbreak}

\newcommand{\ignore}[1]{}
\usepackage[numbers,sort&compress]{natbib}

\begin{document}
\baselineskip 16pt

\title{Completing the spectrum of almost resolvable cycle systems with odd cycle length}

\author{\small  L. Wang and H. Cao \thanks{Research
supported by the National Natural Science Foundation of China
under Grant 11571179, and the Priority Academic
Program Development of Jiangsu Higher Education Institutions.
E-mail: {\sf caohaitao@njnu.edu.cn}} \\
\small Institute of Mathematics, \\ \small  Nanjing Normal
University, Nanjing 210023, China}

\date{}
\maketitle
\begin{abstract}

In this paper, we construct
almost resolvable cycle systems of order $4k+1$ for odd $k\ge 11$.
This completes the proof of the existence of
almost resolvable cycle systems with odd cycle length. As a by-product, some new solutions to the Hamilton-Waterloo problem are also obtained.


\smallskip
\end{abstract}

\section{Introduction}

In this paper, we use $V(H)$ and $E(H)$ to denote the vertex-set and
the edge-set of a graph $H$, respectively.   We  denote the
cycle of length $k$ by $C_k$ and the complete graph on $v$ vertices
by $K_v$.
A {\it factor} of a graph $H$ is a spanning
subgraph whose vertex-set coincides with $V(H)$.
If its connected components are isomorphic to $G$, we call it a {\it
$G$-factor}.
A {\it $G$-factorization} of $H$ is a set of edge-disjoint $G$-factors of $H$
whose edge-sets partition $E(H)$.
 A $C_k$-factorization of $H$ is a partition of
$E(H)$ into $C_k$-factors.
An $r$-regular factor is called an {\it $r$-factor}.
Also, a {\it $2$-factorization} of a graph $H$ is a partition of $E(H)$
into 2-factors.



A $k$-{\it cycle system} of order $v$ is  a collection of
$k$-cycles which partition $E(K_v)$. A $k$-cycle system
of order $v$ exists if and only if $ 3 \leq k \leq v$, $v \equiv 1
\pmod 2$ and $v(v-1) \equiv 0 \pmod {2k}$ \cite{AG,SM}.
A $k$-cycle system of order $v$ is {\it resolvable} if it has a $C_k$-factorization.
A resolvable $k$-cycle system of order $v$
exists if and only if $3 \leq k \leq v$, $v$ and $k$ are odd, and
$v \equiv 0 \pmod k$, see \cite{AH,ASSW, HS, L, L2,PWL, R}. If $v \equiv 1 \pmod {2k}$, then a
$k$-cycle system exists, but it is not resolvable. In this case, Vanstone et al.
\cite{VSS} started the research of the existence of an almost
resolvable $k$-cycle system.

In a $k$-cycle system of order $v$, a collection of $(v-1)/k$ disjoint $k$-cycles is called an {\it
almost parallel class}. In a $k$-cycle system of order $v \equiv 1
\pmod {2k}$, the maximum possible number of almost parallel
classes is $(v-1)/2$, in which case a  half-parallel class
 containing $(v-1)/2k$ disjoint $k$-cycles is left over. A
$k$-cycle system of order $v$ whose cycle set can be partitioned
into $(v-1)/2$ almost parallel classes and a half-parallel class
is called an {\it almost resolvable $k$-cycle system}, denoted by
$k$-ARCS$(v)$.

For recursive constructions of almost resolvable $k$-cycle
systems, C. C. Lindner, et al. \cite{LMR} have considered the
general existence problem of almost resolvable $k$-cycle system
from the commutative quasigroup for $k\equiv0\pmod2$ and made a
 hypothesis:
 if there exists a {\rm $k$-ARCS$(2k+1)$} for $k\equiv0\pmod2$
 and $k\geq8$, then  there exists a {\rm
$k$-ARCS$(2kt+1)$} except possibly for $t=2$. H. Cao et al.
\cite{CNT,NC,WLC} continued to consider the recursive constructions of
an almost resolvable $k$-cycle system  for $k\equiv1\pmod2$.
Many authors contributed to the following known results.


\begin{theorem} {\rm (\cite{ABHL,BHL,CNT,DLM,DLR,LMR,VSS})}
\label{3-14} Let $k \geq 3$, $t \geq 1$ be integers and $n = 2kt +
1$. There exists a $k$-{\rm ARCS}$(n)$ for $k \in \{3, 4, 5, 6, 7,
8, 9, 10, 14\}$, except for $(k, n) \in \{(3, 7), (3, 13), (4,
9)\}$ and except possibly for $(k, n) \in \{(8, 33), (14, 57)\}$.
\end{theorem}

\begin{theorem} {\rm (\cite{NC,WLC})}
\label{2kt+1} For any odd $k\geq 11$, there exists a {\rm
$k$-ARCS$(2kt+1)$}, where $t\ge 1$ and $t\neq 2$.
\end{theorem}

In this paper, we construct
almost resolvable cycle systems of order $4k+1$ for odd $k\ge 11$.
Combining the known results in Theorems~\ref{3-14}-\ref{2kt+1}, we
will prove the following main result.

\begin{theorem}
\label{main} For any odd $k\geq 3$, there exists a {\rm
$k$-ARCS$(2kt+1)$} for  all $t\ge 1$ except for $(k, t) \in \{(3, 1), (3, 2)\}$.
\end{theorem}

\section{Preliminary}

In this section we present a basic lemma for the construction of a $k$-ARCS$(4k+1)$. The main idea is to find some
initial cycles with special properties such that all the required
almost parallel classes can be obtained from them. We  need
the following notions for that lemma.

Suppose $\Gamma$ is an additive group and  $I=\{\infty_1,\infty_2,\ldots,\infty_f\}$ is a set
which is disjoint with $\Gamma$.
We will consider an action of
$\Gamma$ on $\Gamma \ \cup \ I$ which coincides with the
right regular action on the elements of $\Gamma$, and the action of $\Gamma$ on $I$ will coincide with the identity map.
In other words, for any $\gamma\in \Gamma$, we have that $x+\gamma$ is the image under $\gamma$ of any
$x \in \Gamma$, and $x+\gamma=x$ holds for any $x \in I$.
Given a graph $H$ with vertices in $\Gamma\ \cup\ I$, the \emph{translate}
of $H$ by an element $\gamma$ of $\Gamma$ is the graph $H+\gamma$
obtained from $H$ by replacing each vertex $x\in V(H)$ with the
vertex $x+\gamma$. The {\it development} of $H$ under a subgroup
$\Sigma$ of  $\Gamma$ is the collection
$dev_\Sigma(H)=\{H+x\;|\;x\in \Sigma\}$ of all translates of $H$
by an element of $\Sigma$.


For our constructions, we set  $\Gamma=\mathbb{Z}_{k}\times \mathbb{Z}_{4}$.
Given a graph $H$ with vertices in $\Gamma$ and any pair $(r,s)\in \mathbb{Z}_{4}\times
\mathbb{Z}_{4}$, we set $\Delta_{(r,s)}H=\{x-y \ | \  \{(x,r),(y,s) \} \in E(H)\}$.
Finally, given a list $\Hh = \{H_1,H_2,\ldots,H_t\}$ of graphs, we denote by
$\Delta_{(r,s)} \Hh=\cup_{i=1}^t \Delta_{(r,s)} H_i$ the multiset union of the $\Delta_{(r,s)} H_i$s.

\begin{lemma}\label{A}
Let $v=4k+1$ and $\C =\{F_1, F_2\}$ where each $F_i\ (i=1,2)$ is a vertex-disjoint union of four cycles of length $k$ satisfying the following conditions:\\
$(i)$ $V(F_i)=((\mathbb{Z}_{k}\times \mathbb{Z}_{4})\cup\{\infty\})\backslash \{(a_i,b_i)\}$ for some $(a_i,b_i)\in \mathbb{Z}_{k}\times \mathbb{Z}_{4}$, $i=1,2$;\\
$(ii)$ $\infty$  has a neighbor in $\mathbb{Z}_{k}\times\{j\}$ for each $j \in \mathbb{Z}_{4}$;\\
$(iii)$ $\Delta_{(p, p)}\C = \mathbb{Z}_{k}\setminus\{0\}$ for each $p \in \{0,1\}$;\\
$(iv)$ $\Delta_{(q, q)}\C = \mathbb{Z}_{k}\setminus\{0,\pm d_q\}$ for each $q \in \{2,3\}$, where $d_q$ satisfies $(d_q,k)=1$;\\
$(v)$ $\Delta_{(r, s)}\C = \mathbb{Z}_{k}$ for each pair $(r,s)\in \mathbb{Z}_{4}\times
\mathbb{Z}_{4}$ satisfying $r\not= s$.\\
Then, there exists a $k$-{\rm ARCS}$(v)$.
\end{lemma}

\noindent{\it Proof:\ } Let $V(K_{v})=(\mathbb{Z}_{k}\times
\mathbb{Z}_{4})\cup\{\infty\}$.  Note that $0,d_q,2d_q,\ldots,(k-1)d_q$ are $k$ distinct elements since $(d_q,k)=1$.
Then we have the required half parallel class which is formed by the two cycles
$((0,q),(d_q,q)$, $(2d_q,q),\ldots,((k-1)d_q,q))$, $q=2,3$. By $(i)$, we know that $F_i$ is an
almost parallel class. All the required $2k$ almost parallel
classes are $F_i+(l,0)$, $i=1,2$, $l \in \mathbb{Z}_k$.

Now we show that the half parallel class and the $2k$ almost
parallel classes form a $k$-ARCS$(v)$.
 Let  $F'$ be a graph with the edge-set $\{\{(a,q),(a+d_q,q)\} \ | \ a\in \mathbb{Z}_k, q=2,3 \}$ and
 $\Sigma:=\mathbb{Z}_{k}\times\{0\}$. Let $\mathcal{F} = dev_{\Sigma}(\C)\ \cup \ F'$.
The total number of edges -- counted with their respective
multiplicities -- covered by the almost parallel classes and the half parallel class of $\cal F$ is
$2k(4k+1)$, that is exactly the size of $E(K_v)$. Therefore, we only
need to prove that every pair of vertices  lies in a suitable
translate of $\C$ or in $F'$. By $(ii)$,  an edge $\{(z, j),
\infty\}$ of ${K}_v$ must appear in a cycle of $ dev_{\Sigma}(\C)$.

Now consider an edge $\{(z, j), (z', j')\}$ of ${K}_v$ whose
vertices both belong to $\mathbb{Z}_{k}\times \mathbb{Z}_{4}$. If $j=j'\in\{2,3\}$ and $z-z'
\in \{ \pm d_q\}$, then this edge belongs to $F'$. In all other cases there
is,  by $(iii)$-$(v)$, an edge of some $F_i$ of the form $\{(w, j), (w', j')\}$
such that $w-w'=z-z'$. It then follows that $F_i+(-w'+z', 0)$ is an
almost parallel class of $ dev_{\Sigma}(F_i)$ containing the edge
$\{(z, j), (z', j')\}$ and the conclusion follows.\qed

\section{ $k$-{\rm ARCS}$(4k+1)$ for $k \equiv 1 \pmod 4$}

In this section, we will  prove the existence of a $k$-{\rm ARCS}$(4k+1)$ for $k \equiv 1 \pmod 4$.

\begin{lemma}
\label{4k+1}  For any $ k\geq 13$ and $k \equiv 1 \pmod 4$, there exists a $k$-{\rm ARCS}$(4k+1)$.
\end{lemma}

\noindent {\it Proof:} Let $v=4k+1$ and $k=4n+1$, $n\ge 3$. We use Lemma~\ref{A} to construct a
$k$-{\rm ARCS}$(v)$ with $V(K_{v})=(\mathbb{Z}_{k}\times
\mathbb{Z}_{4})\cup\{\infty\}$.
The required parameters in $(i)$ and $(iv)$ of Lemma~\ref{A} are $(a_1,b_1)=(0,3)$, $(a_2,b_2)=(0,2)$, $d_2=2$, and $d_3=\frac{k-1}{2}$.
The required 8 cycles in $F_1=\{C_1, C_2, C_3, C_4\}$ and $F_2=\{C_5, C_6, C_7, C_8\}$ are listed as below.

The cycle $C_1$ is the concatenation of the sequences
$T_1$, $(0,0)$, and $T_2$, where\\
{\footnotesize
$T_1=((n,0),(-n,1),\ldots,\underline{(n-i,0),(-(n-i),1)},\ldots,(1,0),(-1,1))$, $0 \leq i \leq n-1$;\\
$T_2=((1,1),(-1,0),\ldots,\underline{(1+i,1),(-(1+i),0)},\ldots,(n,1),(-n,0))$, $0 \leq i \leq n-1$.
}

\vspace{5pt}

\noindent {\bf Note:} Actually $T_1$ can be viewed as  the concatenation of the sequences
$T_1^0, T_1^1,\dots, T_1^{n-1}$, where the general formula is $T_1^i=((n-i,0),(-(n-i),1))$, $0 \leq i \leq n-1$.
Thus, for brevity, we just list the first sequence at the beginning of $T_1$ and the last sequence at
 the end of $T_1$, and use the underlined sections to give the general formula in the middle of $T_1$.
 We give the range of $i$ after the sequence $T_1$ so that the reader can easily calculate the number of vertices in $T_1$.
Similarly, this partial underlining happens ahead in some places as well.

\vspace{5pt}

The cycle $C_2$ is the concatenation of the sequences
$T_1$, $T_2$, and $(0,2)$, where \\
{\footnotesize
$T_1=((1,2),(-1,3),\ldots,\underline{(1+i,2),(-(1+i),3)},\ldots,(n,2),(-n,3))$, $0 \leq i \leq n-1$;\\
$T_2=((-n,2),(n,3),\ldots,\underline{(-(n-i),2),(n-i,3)},\ldots,(-1,2),(1,3))$, $0 \leq i \leq n-1$.
}

The cycle $C_3$ is the concatenation of the sequences
$T_1$, $T_2$, and $T_3$, where \\
{\footnotesize
$T_1=((n+1,1),(-(n+1),2),\ldots,\underline{(n+1+i,1),(-(n+1+i),2)},\ldots,(2n-1,1),(-(2n-1),2))$, $0 \leq i \leq n-2$;\\
$T_2=((-2n,1),(2n,2),\ldots,\underline{(-(2n-i),1),(2n-i,2)},\ldots,(-(n+1),1),(n+1,2))$, $0 \leq i \leq n-1$;\\
$T_3=((-2n,0),(-2n,2),(0,1))$.
}

The cycle $C_4$ is the concatenation of the sequences $\infty$,
$T_1$, $T_2$, and $T_3$ (By a slight abuse of notation, here $\infty$ is regarded as a sequence.), where \\
{\footnotesize
$T_1=((-2n,3),(2n,3),\ldots,\underline{(-(2n-i),3),(2n-i,3)},\ldots,(-(n+1),3),(n+1,3))$, $0 \leq i \leq n-1$;\\
$T_2=((n+1,0),(-(n+1),0),\ldots,\underline{(n+1+i,0),(-(n+1+i),0)},\ldots,(2n-1,0),(-(2n-1),0))$, $0 \leq i \leq n-2$;\\
$T_3=((2n,0),(2n,1))$.
}

The cycle $C_5$ is the concatenation of the sequences
$T_1$, $(0,1)$, and $T_2$, where \\
{\footnotesize
$T_1=((n,1),(-n,3),\ldots,\underline{(n-i,1),(-(n-i),3)},\ldots,(1,1),(-1,3))$, $0 \leq i \leq n-1$;\\
$T_2=((1,3),(-1,1),\ldots,\underline{(1+i,3),(-(1+i),1)},\ldots,(n,3),(-n,1))$, $0 \leq i \leq n-1$.
}

The cycle $C_6$ is the concatenation of the sequences
$T_1$, $T_2$, $T_3$, and $(-n,0)$, where \\
{\footnotesize
$T_1=((-(n+1),0),(n+1,3),\ldots,\underline{(-(n+1+i),0),(n+1+i,3)},\ldots,(-(2n-1),0),(2n-1,3))$, $0 \leq i \leq n-2$;\\
$T_2=((-2n,0),(0,3))$;\\
$T_3=((2n,0),(-2n,3),\ldots,\underline{(2n-i,0),(-(2n-i),3)},\ldots,(n+1,0),(-(n+1),3))$, $0 \leq i \leq n-1$.
}

Next,  we consider the cycles $C_7$ and $C_8$.

For $k=13,25$, the two cycles are listed as follows. \\
{\footnotesize
$k=13$:\\
$C_7=( ( 1, 2),(-5, 1),(-5, 2),( 0, 0),( 3, 2),(-2, 0),( 4, 2),( 3, 0),( 5, 2),( 1, 0),(-3, 2),(-1, 0),(-2, 2) )$;\\
$C_8=(\infty, ( 2, 0),(-4, 2),(-6, 1),( 5, 1),( 4, 1),(-4, 1),( 6, 1),( 6, 3),(-6, 2),( 2, 2),( 6, 2),(-1, 2))$.
}\\
{\footnotesize
$k=25$:\\
$C_7=(( 1, 2),(-11, 1),(-11, 2),( 0, 0),( 2, 2),( 1, 0),( 4, 2),( 5, 0),( 3, 2),( 6, 0),(10, 2),( 2, 0),( 7, 2),(-5, 0),( 5, 2),$ \\
\hspace*{0.9cm} $(-4, 0),(11, 2),( 4, 0),(-8, 2),(-1, 0),(-10, 2),(-2, 0),(-7, 2),(-3, 0),( 8, 2))$;\\
$C_8=(\infty, ( 3, 0),( 9, 2),( 7, 1),(-7, 1),( 8, 1),(-8, 1),( 9, 1),(11, 1),(-10, 1),(10, 1),(-9, 1),(-12, 1),(12, 1),$ \\
\hspace*{0.9cm} $(12, 3),(-12, 2),(-6, 2),( 6, 2),(-3, 2),(12, 2),(-2, 2),(-5, 2),(-1, 2),(-9, 2),(-4, 2))$.
}

For $ k\geq 17$ and $k \neq 25$, we distinguish the following two cases.

\vspace{5pt}

{\bf Case 1:} $k\equiv5\pmod{8}$ and $ k\geq 21$.

In this case, $C_7$ is the concatenation of the sequences
$S_1$, $S_2$, $\ldots$, $S_6$ as follows.\\
{\footnotesize
$S_1=((-(2n-3),2),(-(2n-1),1))$;\\
$S_2=((1,2),(-1,0),\ldots,\underline{(1+i,2),(-(1+i),0)},\ldots,(\frac{n-1}{2},2),(-\frac{n-1}{2},0))$, $0 \leq i \leq \frac{n-3}{2}$;\\
$S_3=((-\frac{n+1}{2},2),(\frac{n+1}{2},0),\ldots,\underline{(-(\frac{n+1}{2}+i),2),(\frac{n+1}{2}+i,0)},\ldots,(-(n-1),2),(n-1,0))$, $0 \leq i \leq \frac{n-3}{2}$;\\
$S_4=((2n-1,2),(0,0), (-(2n-1),2))$;\\
$S_5=((1,0),(-1,2),\ldots,\underline{(1+i,0),(-(1+i),2)},\ldots,(\frac{n-1}{2},0),(-\frac{n-1}{2},2))$, $0 \leq i \leq \frac{n-3}{2}$;\\
$S_6=((-\frac{n+1}{2},0),(\frac{n+1}{2},2),\ldots,\underline{(-(\frac{n+1}{2}+i),0),(\frac{n+1}{2}+i,2)},\ldots,(-(n-1),0),(n-1,2))$, $0 \leq i \leq \frac{n-3}{2}$.
}

For the last cycle $C_8$, when $k=21, 29$, it is listed as below respectively.\\
{\footnotesize
$k=21$:\\
$C_8=(\infty, ( 5, 0),(-6, 2),(-6, 1),( 6, 1),(-7, 1),( 7, 1),( 9, 1),(-8, 1),( 8, 1),(-10, 1),(10, 1),(10, 3),(-10, 2),$ \\
\hspace*{0.9cm} $( 8, 2),(-5, 2),( 7, 2),(-8, 2),( 6, 2),(10, 2),( 5, 2)
)$.
}\\
{\footnotesize
$k=29$:\\
$C_8=(\infty, ( 7, 0),(-8, 2),(-8, 1),( 8, 1),(-9, 1),( 9, 1),(-10, 1),(10, 1),(-14, 1),(-12, 1),(13, 1),(12, 1),(-11, 1),$ \\
\hspace*{0.9cm} $(11, 1),(14, 1),(14, 3),(-14, 2),( 7, 2),(11, 2),(14, 2),( 8, 2),(-7, 2),( 9, 2),(-9, 2),(10, 2),(-10, 2),(12, 2),$
\hspace*{0.9cm} $(-12, 2))$.
}

When $k\geq 37$, $C_8$ is the concatenation of the sequences $\infty$,
$T_1$, $T_2$, $\ldots$, $T_8$, where \\
{\footnotesize
$T_1=((n,0),(-(n+1),2))$;\\
$T_2=((-(n+1),1),(n+1,1),\ldots,\underline{(-(n+1+i),1),(n+1+i,1)},\ldots,(-\frac{3n-1}{2},1),(\frac{3n-1}{2},1))$, $0 \leq i \leq \frac{n-3}{2}$.
}

For the other 6 sequences $T_3, T_4,\ldots, T_8$, we distinguish the following 3 subcases.

{\bf Case 1.1:} $k\equiv5\pmod{24}$ and $k\geq 53$.\\
{\footnotesize
$T_3=((\frac{3n+5}{2},1),(-\frac{3n+5}{2},1),(\frac{3n+3}{2},1),(-\frac{3n+3}{2}),1),(\frac{3n+1}{2},1),(-\frac{3n+1}{2},1),\ldots, \underline{(\frac{3n+5}{2}+3i,1),(-(\frac{3n+5}{2}+3i),1),}$
\hspace*{1cm} $\underline{(\frac{3n+3}{2}+3i,1),(-(\frac{3n+3}{2}+3i),1),(\frac{3n+1}{2}+3i,1),(-(\frac{3n+1}{2}+3i),1)},\ldots,(2n-7,1),(-(2n-7),1),$\\
\hspace*{1cm}  $(2n-8,1),(-(2n-8),1),(2n-9,1),(-(2n-9),1))$, $0 \leq i \leq \frac{n-19}{6}$;\\
$T_4=((2n-4,1),(-(2n-3),1),(2n-3,1),(-(2n-2),1),(2n-2,1),(-(2n-6),1),(2n-6,1),(-(2n-5),1)$,
\hspace*{1cm} $(2n-5,1),(-(2n-4),1),(-2n,1),(2n-1,1),(2n,1),(2n,3),(-2n,2),(n,2),(-n,2))$;\\
$T_5=((n+2,2),(-(n+2),2),(n+1,2),(-(n+4),2),(n+3,2),(-(n+3),2),$ $\ldots$, $\underline{(n+2+3i,2),(-(n+2+3i),2),}$
 \hspace*{1cm} $\underline{(n+1+3i,2),(-(n+4+3i),2),(n+3+3i,2),(-(n+3+3i),2)},$ $\ldots,(\frac{3n-9}{2},2),(-\frac{3n-9}{2},2),(\frac{3n-11}{2},2),$ \\
 \hspace*{1cm} $(-\frac{3n-5}{2},2),(\frac{3n-7}{2},2),(-\frac{3n-7}{2},2))$, $0 \leq i \leq \frac{n-13}{6}$;\\
$T_6=((\frac{3n-3}{2},2),(-\frac{3n-3}{2},2),(\frac{3n-1}{2},2),(-\frac{3n-1}{2},2),(\frac{3n+3}{2},2),$ $(\frac{3n-5}{2},2))$;\\
$T_7=((\frac{3n+5}{2},2),(-\frac{3n+3}{2},2),(\frac{3n+1}{2},2),(-\frac{3n+5}{2},2),(\frac{3n+9}{2},2),$ $(-\frac{3n+1}{2},2),
\ldots,\underline{(\frac{3n+5}{2}+3i,2),(-(\frac{3n+3}{2}+3i),2),}$
 \hspace*{1cm} $\underline{(\frac{3n+1}{2}+3i,2),(-(\frac{3n+5}{2}+3i),2),(\frac{3n+9}{2}+3i,2),(-(\frac{3n+1}{2}+3i),2)},\ldots, (2n-7,2),(-(2n-8),2),$\\
 \hspace*{1cm} $(2n-9,2),(-(2n-7),2),(2n-5,2),(-(2n-9),2))$, $0 \leq i \leq \frac{n-19}{6}$;\\
$T_8=((2n-4,2),(-(2n-2),2),(2n,2),(-(2n-5),2),(2n-6,2),(-(2n-4),2),(2n-3,2),(-(2n-6),2)$, \hspace*{1cm} $(2n-2,2))$.
}

{\bf Case 1.2:} $k\equiv13\pmod{24}$ and $ k\geq 37$.\\
{\footnotesize
$T_3=((\frac{3n+5}{2},1),(-\frac{3n+5}{2},1),(\frac{3n+3}{2},1),(-\frac{3n+3}{2},1),(\frac{3n+1}{2},1),(-\frac{3n+1}{2},1),\ldots, \underline{(\frac{3n+5}{2}+3i,1),(-(\frac{3n+5}{2}+3i),1),}$
\hspace*{1cm} $\underline{(\frac{3n+3}{2}+3i,1),(-(\frac{3n+3}{2}+3i),1),(\frac{3n+1}{2}+3i,1),(-(\frac{3n+1}{2}+3i),1)},\ldots,(2n-2,1),(-(2n-2),1),$\\
\hspace*{1cm} $(2n-3,1),(-(2n-3),1),(2n-4,1),(-(2n-4),1))$, $0 \leq i \leq \frac{n-9}{6}$;\\
$T_4=((-2n,1),(2n-1,1),(2n,1),(2n,3),(-2n,2),(n,2),(-n,2))$;\\
$T_5=((n+2,2),(-(n+2),2),(n+1,2),(-(n+4),2),(n+3,2),(-(n+3),2),$ $\ldots$, $\underline{(n+2+3i,2),(-(n+2+3i),2),}$
 \hspace*{1cm} $ \underline{(n+1+3i,2),(-(n+4+3i),2),(n+3+3i,2),(-(n+3+3i),2)},$ $\ldots,(\frac{3n-11}{2},2),(-\frac{3n-11}{2},2),$\\
 \hspace*{1cm} $(\frac{3n-13}{2},2),(-\frac{3n-7}{2},2),(\frac{3n-9}{2},2),(-\frac{3n-9}{2},2))$, $0 \leq i \leq \frac{n-15}{6}$;\\
$T_6=((\frac{3n-5}{2},2),(-\frac{3n-5}{2},2),(\frac{3n-7}{2},2),(\frac{3n+5}{2},2),(\frac{3n-3}{2},2)$, $(-\frac{3n-3}{2},2),(\frac{3n-1}{2},2))$;\\
$T_7=((-\frac{3n-1}{2},2),(\frac{3n+11}{2},2),(-\frac{3n+1}{2},2),(\frac{3n+1}{2},2),$ $(-\frac{3n+3}{2},2),(\frac{3n+3}{2},2),\ldots, \underline{(-(\frac{3n-1}{2}+3i),2),(\frac{3n+11}{2}+3i,2),}$\\
 \hspace*{1cm} $\underline{(-(\frac{3n+1}{2}+3i),2),(\frac{3n+1}{2}+3i,2),(-(\frac{3n+3}{2}+3i),2),(\frac{3n+3}{2}+3i,2)},\ldots,(-(2n-8),2),(2n-2,2),$\\
 \hspace*{1cm} $(-(2n-7),2),(2n-7,2)$, $(-(2n-6),2)),(2n-6,2)$, $0 \leq i \leq \frac{n-15}{6}$;\\
$T_8=((-(2n-5),2),(2n-3,2),(-(2n-4),2),(2n,2),(-(2n-2),2),(2n-4,2))$.
}

{\bf Case 1.3:} $k\equiv21\pmod{24}$ and  $ k\geq 45$.\\
{\footnotesize
$T_3=((\frac{3n+5}{2},1),(-\frac{3n+5}{2},1),(\frac{3n+3}{2},1),(-\frac{3n+3}{2},1),(\frac{3n+1}{2},1),(-\frac{3n+1}{2},1),\ldots, \underline{(\frac{3n+5}{2}+3i,1),(-(\frac{3n+5}{2}+3i),1),}$\\
\hspace*{1cm} $\underline{(\frac{3n+3}{2}+3i,1),(-(\frac{3n+3}{2}+3i),1),(\frac{3n+1}{2}+3i,1),(-(\frac{3n+1}{2}+3i),1)},\ldots,(2n-6,1),(-(2n-6),1),$\\
\hspace*{1cm} $(2n-7,1),(-(2n-7),1),(2n-8,1),(-(2n-8),1))$, $0 \leq i \leq \frac{n-17}{6}$;\\
$T_4=((2n-3,1),(-(2n-3),1),(2n-2,1),(-(2n-2),1),(2n-5,1),(-(2n-5),1),(2n-4,1),(-(2n-4),1)$,
\hspace*{1cm} $(-2n,1),(2n-1,1),(2n,1),(2n,3),(-2n,2),(n,2),(-n,2))$;\\
$T_5=((n+2,2),(-(n+2),2),(n+1,2),(-(n+4),2),(n+3,2),(-(n+3),2),$ $\ldots,$ $\underline{(n+2+3i,2),(-(n+2+3i),2),}$
 \hspace*{1cm} $ \underline{(n+1+3i,2),(-(n+4+3i),2),(n+3+3i,2),(-(n+3+3i),2)},\ldots,(\frac{3n-13}{2},2),(-\frac{3n-13}{2},2),$\\
 \hspace*{1cm} $(\frac{3n-15}{2},2),(-\frac{3n-9}{2},2),(\frac{3n-11}{2},2),(-\frac{3n-11}{2},2))$, $0 \leq i \leq \frac{n-17}{6}$;\\
$T_6=((\frac{3n-7}{2},2),(-\frac{3n-7}{2},2),(\frac{3n-9}{2},2),(-\frac{3n-3}{2},2),(\frac{3n-1}{2},2),(-\frac{3n-5}{2},2),(\frac{3n-5}{2},2),$
$(\frac{3n+3}{2},2),(\frac{3n-3}{2},2),$\\
 \hspace*{1cm} $(-\frac{3n+1}{2},2),(\frac{3n+1}{2},2))$;\\
$T_7=((-\frac{3n+3}{2},2),(\frac{3n+7}{2},2),(-\frac{3n-1}{2},2),(\frac{3n+9}{2},2),$ $(-\frac{3n+7}{2},2),(\frac{3n+5}{2},2),\ldots, $
$\underline{(-(\frac{3n+3}{2}+3i),2),(\frac{3n+7}{2}+3i,2),}$
 \hspace*{1cm} $\underline{(-(\frac{3n-1}{2}+3i),2),(\frac{3n+9}{2}+3i,2),(-(\frac{3n+7}{2}+3i),2),(\frac{3n+5}{2}+3i,2)},\ldots,(-(2n-7),2),(2n-5,2),$\\
 \hspace*{1cm} $(-(2n-9),2),(2n-4,2),(-(2n-5),2)$, $(2n-6,2))$, $0 \leq i \leq \frac{n-17}{6}$;\\
$T_8=((-(2n-4),2),(2n-3,2),(-(2n-2),2),(2n-2,2),(-(2n-6),2),(2n,2))$.
}

{\bf Case 2:} $k\equiv1\pmod{8}$, $ k\geq 17$ and $k\neq 25$.

In this case, $C_7$ is  the concatenation of the sequences
$S_1$, $S_2$, $(-n,2)$, $S_3$, $S_4$, $S_5$, and $S_6$.\\
{\footnotesize
$S_1=((-(2n-3),2),(-(2n-1),1))$;\\
$S_2=((1,2),(-1,0),\ldots,\underline{(1+i,2),(-(1+i),0)},\ldots,(n-1,2),(-(n-1),0))$, $0 \leq i \leq n-2$;\\
$S_3=((1,0),(-1,2),\ldots,\underline{(1+i,0),(-(1+i),2)},\ldots,(\frac{n-2}{2},0),(-\frac{n-2}{2},2))$, $0 \leq i \leq \frac{n-4}{2}$;\\
$S_4=((\frac{n}{2},0),(-\frac{3n}{2},2))$;\\
$S_5=((\frac{n+2}{2},0),(-\frac{n+2}{2},2),\ldots,\underline{(\frac{n+2}{2}+i,0),(-(\frac{n+2}{2}+i),2)},\ldots,(n-1,0),(-(n-1),2))$, $0 \leq i \leq \frac{n-4}{2}$;\\
$S_6=((n,0),(n+1,2))$.
}

For the cycle $C_8$, when $k=17,33,41$, it is listed as below respectively.\\
{\footnotesize
$k=17$:\\
$C_8=(\infty, ( 0, 0),( 7, 2),( 7, 1),( 6, 1),(-8, 1),(-6, 1),( 5, 1),(-5, 1),( 8, 1),( 8, 3),(-8, 2),( 6, 2),(-2, 2),(-7, 2),$ \\
\hspace*{0.9cm} $( 4, 2),( 8, 2))$.
}\\
{\footnotesize
$k=33$:\\
$C_8=(\infty, ( 0, 0),(15, 2),(15, 1),( 9, 1),(-9, 1),(10, 1),(-10, 1),(11, 1),(-11, 1),(12, 1),(13, 1),(-12, 1),(14, 1),$ \\
\hspace*{0.9cm} $(-14, 1),(-16, 1),(-13, 1),(16, 1),(16, 3),(-16, 2),( 8, 2),(12, 2),(-14, 2),(16, 2),(-11, 2),(14, 2),$ \\
\hspace*{0.9cm} $(-9, 2),(10, 2),(-10, 2),(11, 2),(-4, 2),(13, 2),(-15, 2))$.
}\\
{\footnotesize
$k=41$:\\
$C_8=(\infty,( 0, 0),(19, 2),(19, 1),(-12, 1),(13, 1),(-13, 1),(11, 1),(-11, 1),(12, 1),(-15, 1),(14, 1),(-14, 1),$\\
\hspace*{1cm} $(18,1),(15, 1),(16, 1),(-17, 1),(17, 1),(-18, 1),(-20, 1),(-16, 1),(20, 1),(20, 3),(-20, 2),(10, 2),$\\
\hspace*{1cm} $(13, 2),(-13, 2),(14, 2),(-11, 2),(12, 2),(-12, 2),(17, 2),(-5, 2),(16, 2),(-19, 2),(18, 2),(-14, 2),$\\
\hspace*{1cm} $ (20, 2),(-16,2),(15, 2),(-18, 2)).$
}

When $k\geq 49$, the cycle $C_8$ is the concatenation of the sequences $\infty$,
$T_1$, $T_2$, $\ldots$, $T_9$, where \\
{\footnotesize
$T_1=((0,0),(2n-1,2),(2n-1,1))$.}\\
For the other 8 sequences $T_2$, $T_3$, $\ldots$, $T_9$, we distinguish 3 subcases.

{\bf Case 2.1:} $k\equiv 1\pmod{24}$ and  $ k\geq 49$.\\
{\footnotesize
$T_2=((-(n+2),1),(n+3,1),(-(n+3),1),(n+1,1),(-(n+1),1),(n+2,1),\ldots,\underline{(-(n+2+3i),1),(n+3+3i,1),}$
\hspace*{1cm} $ \underline{(-(n+3+3i),1),(n+1+3i,1),(-(n+1+3i),1),(n+2+3i,1)},\ldots,(-\frac{3n-8}{2},1),(\frac{3n-6}{2},1),(-\frac{3n-6}{2},1),$
\hspace*{1cm} $(\frac{3n-10}{2},1),(-\frac{3n-10}{2},1),(\frac{3n-8}{2},1))$, $0 \leq i \leq \frac{n-12}{6}$;\\
$T_3=((-\frac{3n-2}{2},1),(\frac{3n}{2},1),(\frac{3n-2}{2},1),(-\frac{3n-4}{2},1),(\frac{3n-4}{2},1),(-\frac{3n}{2},1),(\frac{3n+4}{2},1),$
$(-\frac{3n+2}{2},1),(-\frac{3n+6}{2},1),(\frac{3n+2}{2},1))$;\\
$T_4=((-\frac{3n+8}{2},1),(\frac{3n+8}{2},1),(-\frac{3n+4}{2},1),(\frac{3n+10}{2},1),(-\frac{3n+12}{2},1),(\frac{3n+6}{2},1),\ldots, \underline{(-(\frac{3n+8}{2}+3i),1),(\frac{3n+8}{2}+3i,1),}$
\hspace*{1cm}  $\underline{(-(\frac{3n+4}{2}+3i),1),(\frac{3n+10}{2}+3i,1),(-(\frac{3n+12}{2}+3i),1),(\frac{3n+6}{2}+3i,1)},\ldots,(-(2n-5),1),(2n-5,1),$\\
\hspace*{1cm}  $(-(2n-7),1),(2n-4,1),(-(2n-3),1),(2n-6,1))$, $0 \leq i \leq \frac{n-18}{6}$;\\
$T_5=((-(2n-2),1),(2n-3,1),(-2n,1),(2n-2,1),(-(2n-4),1),(2n,1),(2n,3),(-2n,2),(n,2))$;\\
$T_6=((n+4,2),(-(n+3),2),(n+3,2),(-(n+2),2),(n+2,2),(-(n+1),2),\ldots,\underline{(n+4+3i,2),(-(n+3+3i),2),}$
 \hspace*{1cm} $\underline{(n+3+3i,2),(-(n+2+3i),2),(n+2+3i,2),(-(n+1+3i),2)},\ldots,(\frac{3n-4}{2},2),(-\frac{3n-6}{2},2),(\frac{3n-6}{2},2),$\\
 \hspace*{1cm} $(-\frac{3n-8}{2},2),(\frac{3n-8}{2},2),(-\frac{3n-10}{2},2))$, $0 \leq i \leq \frac{n-12}{6}$;\\
$T_7=((\frac{3n+2}{2},2),(-\frac{n}{2},2),(\frac{3n-2}{2},2),(-\frac{3n-4}{2},2),(\frac{3n+8}{2},2),(-\frac{3n-2}{2},2),(\frac{3n}{2},2))$;\\
$T_8=((-\frac{3n+2}{2},2),(\frac{3n+14}{2},2),(-\frac{3n+4}{2},2),(\frac{3n+4}{2},2),(-\frac{3n+6}{2},2),(\frac{3n+6}{2},2),\ldots,\underline{(-(\frac{3n+2}{2}+3i),2),(\frac{3n+14}{2}+3i,2),}$
 \hspace*{1cm} $\underline{(-(\frac{3n+4}{2}+3i),2),(\frac{3n+4}{2}+3i,2),(-(\frac{3n+6}{2}+3i),2),(\frac{3n+6}{2}+3i,2)},\ldots,(-(2n-8),2),(2n-2,2),$\\
 \hspace*{1cm} $(-(2n-7),2),(2n-7,2),(-(2n-6),2),(2n-6,2))$, $0 \leq i \leq \frac{n-18}{6}$;\\
$T_9=((-(2n-5),2),(2n,2),(-(2n-2),2),(2n-4,2),(-(2n-4),2),(2n-3,2),(-(2n-1),2))$.
}

{\bf Case 2.2:} $k\equiv9\pmod{24}$ and  $ k\geq 57$.\\
{\footnotesize
$T_2=((-(n+2),1),(n+3,1),(-(n+3),1),(n+1,1),(-(n+1),1),(n+2,1),\ldots, \underline{(-(n+2+3i),1),(n+3+3i,1),}$
\hspace*{1cm} $ \underline{(-(n+3+3i),1),(n+1+3i,1),(-(n+1+3i),1),(n+2+3i,1)},\ldots,(-\frac{3n-4}{2},1),(\frac{3n-2}{2},1),$\\
\hspace*{1cm} $(-\frac{3n-2}{2},1),(\frac{3n-6}{2},1),(-\frac{3n-6}{2},1),(\frac{3n-4}{2},1))$, $0 \leq i \leq \frac{n-8}{6}$;\\
$T_3=((-\frac{3n+2}{2},1),(-\frac{3n}{2},1),(\frac{3n+4}{2},1))$;\\
$T_4=((\frac{3n}{2},1),(-\frac{3n+8}{2},1),(\frac{3n+2}{2},1),(-\frac{3n+4}{2},1),$
$(\frac{3n+10}{2},1),(-\frac{3n+6}{2},1),\ldots,\underline{(\frac{3n}{2}+3i,1),(-(\frac{3n+8}{2}+3i),1),}$\\
\hspace*{1cm}  $ \underline{(\frac{3n+2}{2}+3i,1),(-(\frac{3n+4}{2}+3i),1),(\frac{3n+10}{2}+3i,1),(-(\frac{3n+6}{2}+3i),1)},\ldots,(2n-7,1),(-(2n-3),1),$\\
\hspace*{1cm}  $(2n-6,1),(-(2n-5),1),(2n-2,1),(-(2n-4),1))$, $0 \leq i \leq \frac{n-14}{6}$;\\
$T_5=((2n-4,1),(-2n,1),(2n-3,1),(-(2n-2),1),(2n,1),(2n,3),(-2n,2),(n,2),(n+3,2),(-(n+2),2),$
\hspace*{1cm}  $(n+2,2))$;\\
$T_6=((-(n+1),2),(n+5,2),(-(n+5),2),(n+6,2),$ $(-(n+3),2),(n+4,2),$ $\ldots$, $\underline{(-(n+1+3i),2),(n+5+3i,2),}$
 \hspace*{1cm} $ \underline{(-(n+5+3i),2),(n+6+3i,2),}\underline{(-(n+3+3i),2),(n+4+3i,2)},\ldots,(-\frac{3n-12}{2},2),(\frac{3n-4}{2},2),$ \\
 \hspace*{1cm} $(-\frac{3n-4}{2},2),(\frac{3n-2}{2},2),(-\frac{3n-8}{2},2),(\frac{3n-6}{2},2))$, $0 \leq i \leq \frac{n-14}{6}$;\\
$T_7=((-\frac{3n-6}{2},2),(\frac{3n+4}{2},2),(-\frac{n}{2},2),(\frac{3n}{2},2),(-\frac{3n+2}{2},2))$;\\
$T_8=((\frac{3n+2}{2},2),(-\frac{3n+8}{2},2),(\frac{3n+6}{2},2),(-\frac{3n+6}{2},2),(\frac{3n+10}{2},2),(-\frac{3n-2}{2},2),\ldots,\underline{(\frac{3n+2}{2}+3i,2),(-(\frac{3n+8}{2}+3i),2),}$
 \hspace*{1cm} $\underline{(\frac{3n+6}{2}+3i,2),(-(\frac{3n+6}{2}+3i),2),(\frac{3n+10}{2}+3i,2),(-(\frac{3n-2}{2}+3i),2)},\ldots,(2n-9,2),(-(2n-6),2),$\\
 \hspace*{1cm} $(2n-7,2),(-(2n-7),2),(2n-5,2),(-(2n-11),2))$, $0 \leq i \leq \frac{n-20}{6}$;\\
$T_9=((2n-6,2),(-(2n-4),2),(2n-3,2),(-(2n-8),2),(2n,2),(-(2n-5),2),(2n-4,2),(-(2n-2),2),$
 \hspace*{1cm} $(2n-2,2),(-(2n-1),2))$.
}

{\bf Case 2.3:} $k\equiv 17\pmod{24}$ and $k\geq 65$.\\
{\footnotesize
$T_2=((-(n+2),1),(n+3,1),(-(n+3),1),(n+1,1),(-(n+1),1),(n+2,1),\ldots,\underline{(-(n+2+3i),1),(n+3+3i,1),}$
\hspace*{1cm} $ \underline{(-(n+3+3i),1),(n+1+3i,1),(-(n+1+3i),1),(n+2+3i,1)},\ldots,(-\frac{3n-6}{2},1),(\frac{3n-4}{2},1),$\\
\hspace*{1cm} $(-\frac{3n-4}{2},1),(\frac{3n-8}{2},1),(-\frac{3n-8}{2},1),(\frac{3n-6}{2},1))$, $0 \leq i \leq \frac{n-10}{6}$;\\
$T_3=((-\frac{3n}{2},1),(\frac{3n-2}{2},1),(-\frac{3n-2}{2},1),(\frac{3n+6}{2},1),(\frac{3n}{2},1))$;\\
$T_4=((\frac{3n+2}{2},1),(-\frac{3n+6}{2},1),$
$(\frac{3n+4}{2},1),(-\frac{3n+2}{2},1),(\frac{3n+12}{2},1),(-\frac{3n+4}{2},1),\ldots,\underline{(\frac{3n+2}{2}+3i,1),(-(\frac{3n+6}{2}+3i),1),}$\\
\hspace*{1cm}  $ \underline{(\frac{3n+4}{2}+3i,1),(-(\frac{3n+2}{2}+3i),1),(\frac{3n+12}{2}+3i,1),}\underline{(-(\frac{3n+4}{2}+3i),1)},\ldots,(2n-7,1),(-(2n-5),1),$\\
\hspace*{1cm}  $(2n-6,1),(-(2n-7),1),(2n-2,1),(-(2n-6),1))$, $0 \leq i \leq \frac{n-16}{6}$;\\
$T_5=((2n-4,1),(-(2n-3),1),(2n-3,1),(-(2n-2),1),(-2n,1),(-(2n-4),1),(2n,1),(2n,3),(-2n,2),(n,2),$
\hspace*{1cm}  $(n+3,2),(-(n+3),2),(n+4,2),(-(n+1),2),(n+2,2))$;\\
$T_6=((-(n+2),2),(n+6,2),(-(n+6),2),(n+7,2),(-(n+4),2),(n+5,2),\ldots$, $\underline{(-(n+2+3i),2),(n+6+3i,2),}$
 \hspace*{1cm} $\underline{(-(n+6+3i),2),(n+7+3i,2),(-(n+4+3i),2),(n+5+3i,2)},$ $\ldots,(-\frac{3n-12}{2},2),(\frac{3n-4}{2},2),$\\
 \hspace*{1cm} $(-\frac{3n-4}{2},2),(\frac{3n-2}{2},2),(-\frac{3n-8}{2},2),(\frac{3n-6}{2},2))$, $0 \leq i \leq \frac{n-16}{6}$;\\
$T_7=((-\frac{3n-6}{2},2),(\frac{3n+4}{2},2),(-\frac{n}{2},2),(\frac{3n}{2},2),(-\frac{3n+2}{2},2))$;\\
$T_8=((\frac{3n+2}{2},2),(-\frac{3n+8}{2},2),(\frac{3n+6}{2},2),(-\frac{3n+6}{2},2),(\frac{3n+10}{2},2),(-\frac{3n-2}{2},2),\ldots,\underline{(\frac{3n+2}{2}+3i,2),(-(\frac{3n+8}{2}+3i),2),}$
 \hspace*{1cm} $\underline{(\frac{3n+6}{2}+3i,2),(-(\frac{3n+6}{2}+3i),2),(\frac{3n+10}{2}+3i,2),(-(\frac{3n-2}{2}+3i),2)},\ldots,(2n-10,2),(-(2n-7),2),$\\
 \hspace*{1cm} $(2n-8,2),(-(2n-8),2),(2n-6,2),(-(2n-12),2))$, $0 \leq i \leq \frac{n-22}{6}$;\\
$T_9=((2n-7,2),(-(2n-2),2),(2n-5,2),(-(2n-6),2),(2n-4,2),(-(2n-9),2),(2n-3,2),(-(2n-5),2),$
 \hspace*{1cm} $(2n,2),(-(2n-4),2),(2n-2,2),(-(2n-1),2))$.
}
\qed

\section{ $k$-{\rm ARCS}$(4k+1)$ for $k \equiv 3 \pmod 4$}

In this section, we will prove the existence of a $k$-{\rm ARCS}$(4k+1)$ for $k \equiv 3 \pmod 4$.

\begin{lemma}
\label{4k+3} For any $ k\geq 11$ and $k \equiv 3 \pmod 4$, there exists a $k$-{\rm ARCS}$(4k+1)$.
\end{lemma}

\noindent {\it Proof:} Let $v=4k+1$ and $k=4n+3$, $n\ge 2$. We use Lemma~\ref{A} to construct a
$k$-{\rm ARCS}$(v)$ with $V(K_{v})=(\mathbb{Z}_{k}\times
\mathbb{Z}_{4})\cup\{\infty\}$. Three of the required parameters are $(a_1,b_1)=(0,3)$ and $d_3=\frac{k-1}{2}$.
The other required parameters $a_2,b_2,d_2$ and 8 cycles in $F_1=\{C_1, C_2, C_3, C_4\}$ and $F_2=\{C_5, C_6, C_7, C_8\}$ are listed as below.

The cycle $C_1$ is the concatenation of the sequences
$T_1$, $(0,0)$, and $T_2$, where\\
{\footnotesize
$T_1=((n,0),(-n,1),\ldots,\underline{(n-i,0),(-(n-i),1)},\ldots,(1,0),(-1,1))$, $0 \leq i \leq n-1$;\\
$T_2=((1,1),(-1,0),\ldots,\underline{(1+i,1),(-(1+i),0)},\ldots,(n+1,1),(-(n+1),0))$, $0 \leq i \leq n$.
}

The cycle $C_2$ is the concatenation of the sequences
$T_1$, $T_2$, $T_3$, and $(0,2)$, where \\
{\footnotesize
$T_1=((1,2),(-1,3),\ldots,\underline{(1+i,2),(-(1+i),3)},\ldots,(n,2),(-n,3))$, $0 \leq i \leq n-1$;\\
$T_2=((n+1,2),(n+1,3))$;\\
$T_3=((-n,2),(n,3),\ldots,\underline{(-(n-i),2),(n-i,3)},\ldots,(-1,2),(1,3))$, $0 \leq i \leq n-1$.
}

The cycle $C_3$ is the concatenation of the sequences
$T_1$, $T_2$, and $T_3$, where \\
{\footnotesize
$T_1=((n+2,1),(-(n+2),2),\ldots,\underline{(n+2+i,1),(-(n+2+i),2)},\ldots,(2n,1),(-2n,2))$, $0 \leq i \leq n-2$;\\
$T_2=((-(2n+1),1),(2n+1,2),\ldots,\underline{(-(2n+1-i),1),(2n+1-i,2)},\ldots,(-(n+2),1),(n+2,2))$, $0 \leq i \leq n-1$;\\
$T_3=((-(n+1),1),(-(n+1),2)),(-(2n+1),0),(-(2n+1),2),(0,1))$.
}

The cycle $C_4$ is the concatenation of the sequences $\infty$,
$T_1$, $T_2$, $T_3$, and $T_4$, where \\
{\footnotesize
$T_1=((-(2n+1),3),(2n+1,3),\ldots,\underline{(-(2n+1-i),3),(2n+1-i,3)},\ldots,(-(n+2),3),(n+2,3))$, $0 \leq i \leq n-1$;\\
$T_2=((-(n+1),3),(n+1,0))$;\\
$T_3=((n+2,0),(-(n+2),0),\ldots,\underline{(n+2+i,0),(-(n+2+i),0)},\ldots,(2n,0),(-2n,0))$, $0 \leq i \leq n-2$;\\
$T_4=((2n+1,0),(2n+1,1))$.
}

For $ k=11$, the other four cycles are listed as below.\\
{\footnotesize
$k=11$: $(a_2,b_2)=(\frac{k+1}{2},2)$, $d_2=2$.\\
$C_5=(( 0, 1),( 0, 3),( 1, 1),( 2, 3),( 4, 1),( 1, 3),( 5, 1),(-4, 3),( 2, 1),( 5, 3),(-1, 1))$;\\
$C_6=(( 0, 0),( 3, 3),( 1, 0),(-2, 3),( 3, 0),( 4, 3),( 4, 0),(-3, 3),(-1, 0),(-5, 3),(-4, 0))$;\\
$C_7=(( 2, 2),(-3, 1),( 1, 2),(-5, 0),( 4, 2),(-3, 0),( 5, 2),( 2, 0),( 3, 2),(-2, 0),(-3, 2))$;\\
$C_8=(\infty, ( 5, 0),(-2, 2),(-4, 1),(-2, 1),( 3, 1),(-5, 1),(-1, 3),( 0, 2),(-4, 2),(-1, 2))$.
}

For $k\geq 15$, the two cycles $C_5$ and $C_6$ are defined as below.

The cycle $C_5$ is the concatenation of the sequences
$T_1$, $T_2$, $(0,1)$, $T_3$, and $T_4$, where \\
{\footnotesize
$T_1=((n+2,1),(-(n+1),3))$;\\
$T_2=((n,1),(-n,3),\ldots,\underline{(n-i,1),(-(n-i),3)},\ldots,(1,1),(-1,3))$, $0 \leq i \leq n-1$;\\
$T_3=((1,3),(-1,1),\ldots,\underline{(1+i,3),(-(1+i),1)},\ldots,(n-1,3),(-(n-1),1))$, $0 \leq i \leq n-2$;\\
$T_4=((n,3),(-(n+1),1))$.
}

The cycle $C_6$ is the concatenation of the sequences
$T_1$, $T_2$, $T_3$, and $(n+1,0)$, where \\
{\footnotesize
$T_1=((-(n+2),0),(n+2,3),\ldots,\underline{(-(n+2+i),0),(n+2+i,3)},\ldots,(-(2n+1),0),(2n+1,3))$, $0 \leq i \leq n-1$;\\
$T_2=((2n+1,0),(0,3),(-1,0),(-(2n+1),3))$;\\
$T_3=((2n,0),(-2n,3),\ldots,\underline{(2n-i,0),(-(2n-i),3)},\ldots,(n+2,0),$ $(-(n+2),3))$, $0 \leq i \leq n-2$.
}

For the last two cycles $C_7$ and $C_8$, when $k=15, 19, 23, 27, 35$, they are listed as below.\\
{\footnotesize
$k=15$: $(a_2,b_2)=(\frac{3k-5}{4},2)$, $d_2=4$. \\
$C_7=(( 0, 2),(-7, 1),(-1, 2),( 0, 0),( 1, 2),( 3, 0),( 7, 2),( 1, 0),(-2, 2),(-4, 0),( 3, 2),(-2, 0),( 6, 2),(-3, 0),(-7, 2))$;\\
$C_8=(\infty, ( 2, 0),(-3, 2),(-5, 1),( 6, 1),( 7, 1),(-6, 1),(-3, 1),( 4, 1),( 4, 3),( 5, 2),( 2, 2),(-4, 2),(-6, 2),( 4, 2))$.
}\\
{\footnotesize
$k=19$: $(a_2,b_2)=(\frac{k+1}{2},2)$, $d_2=2$.\\
$C_7=(( 0, 2),(-9, 1),(-7, 2),( 0, 0),( 1, 2),( 2, 0),( 4, 2),( 1, 0),( 7, 2),(-4, 0),(-8, 2),( 3, 0),(-2, 2),( 4, 0),(-6, 2),$
\hspace*{0.9cm} $(-3, 0),( 2, 2),(-5, 0),( 5, 2))$;\\
$C_8=(\infty, (-2, 0),(-4, 2),( 7, 1),(-8, 1),( 9, 1),(-7, 1),(-6, 1),( 8, 1),(-4, 1),( 5, 1),( 5, 3),( 6, 2),( 3, 2),( 9, 2),$\\
\hspace*{0.9cm} $(-3, 2),( 8, 2),(-1, 2),(-5, 2))$.
}\\
{\footnotesize
$k=23$: $(a_2,b_2)=(\frac{3k-1}{4},2)$, $d_2=2$. \\
$C_7=(( 1, 2),(-9, 1),( 2, 2),(-2, 0),( 4, 2),(-3, 0),( 5, 2),(-4, 0),( 6, 2),(-5, 0),( 8, 2),( 5, 0),(-4, 2),( 4, 0),$ \\
\hspace*{0.9cm} $( 3, 2),( 1, 0),(-1, 2),( 3, 0),(-3, 2),( 0, 0),(-5, 2),(-6, 0),(10, 2))$;\\
$C_8=(\infty, ( 2, 0),(-9, 2),(-11, 1),( 8, 1),(-7, 1),(11, 1),(10, 1),(-10, 1),(-8, 1),( 9, 1),(-5, 1),( 6, 1),( 6, 3),$ \\
\hspace*{0.9cm} $( 7, 2),( 0, 2),(-11, 2),(-8, 2),( 9, 2),(-10, 2),(-2, 2),(11, 2),(-7, 2))$.
}\\
{\footnotesize
$k=27$: $(a_2,b_2)=(\frac{k+1}{2},2)$, $d_2=2$. \\
$C_7=(( 0, 2),(-13, 1),(-11, 2),( 0, 0),( 1, 2),( 2, 0),( 4, 2),( 1, 0),( 5, 2),(-7, 0),( 2, 2),( 4, 0),( 9, 2),(-6, 0),( 7, 2),$ \\
\hspace*{0.9cm} $(-4, 0),( 3, 2),(-5, 0),(13, 2),( 3, 0),(-5, 2),(-2, 0),(-8, 2),(-3, 0),(-7, 2),( 6, 0),(-4, 2))$;\\
$C_8=(\infty, ( 5, 0),(-2, 2),(13, 1),(11, 1),(-12, 1),(-11, 1),(-8, 1),( 9, 1),(-9, 1),(12, 1),(-10, 1),(10, 1),(-6, 1),$ \\
\hspace*{0.9cm} $( 7, 1),( 7, 3),( 8, 2),(11, 2),(-6, 2),( 6, 2),(12, 2),(-10, 2),(10, 2),(-3, 2),(-12, 2),(-1, 2),(-9, 2))$.
}\\
{\footnotesize
$k=35$: $(a_2,b_2)=(\frac{k+1}{2},2)$, $d_2=2$. \\
$C_7=(( 1, 2),(-15, 1),( 2, 2),(-2, 0),( 3, 2),(-3, 0),( 4, 2),(-5, 0),( 5, 2),(-6, 0),( 6, 2),(-7, 0),( 7, 2),(-8, 0),$ \\
\hspace*{0.9cm} $( 8, 2),(-9, 0),(11, 2),( 8, 0),( 9, 2),( 7, 0),(-9, 2),( 5, 0),(-8, 2),( 4, 0),(-7, 2),( 3, 0),(-6, 2),(-4, 0),$ \\
\hspace*{0.9cm} $(-5, 2),( 0, 0),(-4, 2),( 2, 0),(-1, 2),( 6, 0),(-2, 2))$;\\
$C_8=(\infty, ( 1, 0),(-16, 2),(17, 1),(-16, 1),(16, 1),(-12, 1),(-17, 1),(-13, 1),(13, 1),(-14, 1),(15, 1),(14, 1),$ \\
\hspace*{0.9cm} $(-10, 1),(11, 1),(-11, 1),(12, 1),(-8, 1),( 9, 1),( 9, 3),(10, 2),(-12, 2),(12, 2),(-11, 2),(14, 2),(-13, 2),$ \\
\hspace*{0.9cm} $(16, 2),(-10, 2),(-15, 2),(13, 2),(-3, 2),(15, 2),( 0, 2),(-14, 2),(17, 2))$.
}

For the last two cycles $C_7$ and $C_8$ of other values,  we distinguish the following 2 cases.

{\bf Case 1:}  $k\equiv 3\pmod{8}$ and $ k\geq 43$. Here, $(a_2,b_2)=(\frac{k+1}{2},2)$, $d_2=2$.

The cycle $C_7$ is the concatenation of the sequences
$S_1$, $S_2$, $\ldots$, $S_5$, $(-\frac{n}{2},0)$, $S_6$, and $S_7$, where\\
{\footnotesize
$S_1=((1,2),(-(2n-1),1))$;\\
$S_2=((2,2),(-2,0),\ldots,\underline{(2+i,2),(-(2+i),0)},\ldots,(\frac{n-2}{2},2),(-\frac{n-2}{2},0))$, $0 \leq i \leq \frac{n-6}{2}$;\\
$S_3=((\frac{n}{2},2),(-\frac{n+2}{2},0),$ $\ldots,\underline{(\frac{n}{2}+i,2),(-(\frac{n+2}{2}+i),0)},\ldots,(n,2),(-(n+1),0))$, $0 \leq i \leq \frac{n}{2}$;\\
$S_4=((n+3,2),(n,0),(n+1,2),(n-1,0),(-(n+1),2))$;\\
$S_5=((n-3,0),(-n,2),\ldots,\underline{(n-3-i,0),(-(n-i),2)},\ldots,$ $(\frac{n-2}{2},0),(-\frac{n+4}{2},2))$, $0 \leq i \leq \frac{n-4}{2}$.
}

For the sequences $S_6$, $S_7$, and the cycle $C_8$, we distinguish 3 subcases.

{\bf Case 1.1:} $k\equiv 3\pmod{24}$ and   $ k\geq 51$.\\
{\footnotesize
$S_6=((-\frac{n+2}{2},2),(\frac{n-8}{2},0),(-\frac{n}{2},2),(\frac{n-4}{2},0),(-\frac{n-8}{2},2),
(\frac{n-6}{2},0),\ldots,\underline{(-(\frac{n+2}{2}-3i),2),(\frac{n-8}{2}-3i,0),(-(\frac{n}{2}-3i),}$\\
\hspace*{1cm}  $\underline{2),(\frac{n-4}{2}-3i,0),(-(\frac{n-8}{2}-3i),2),(\frac{n-6}{2}-3i,0)},\ldots,(-10,2),(5,0),(-9,2),(7,0),(-5,2),(6,0))$, \\
\hspace*{1cm}  $0 \leq i \leq \frac{n-18}{6}$;\\
$S_7=((-7,2),$ $(3,0),(-6,2),(2,0),(-4,2),(0,0),(-3,2),(4,0),(-1,2),(n-2,0),(-2,2))$.
}

The cycle $C_8$ is the concatenation of the sequences $\infty$,
$T_1$, $T_2$, $\ldots$, $T_9$, where\\
{\footnotesize
$T_1=((1,0),(-2n,2),(2n+1,1),(-2n,1),(-(2n+1),1),(-(2n-3),1),(2n-1,1),(-(2n-2),1),(2n,1))$;\\
$T_2=((-(2n-5),1),(2n-2,1),(-(2n-4),1),(2n-4,1),(-(2n-6),1),(2n-3,1),\ldots, \underline{(-(2n-5-3i),1),}$
\hspace*{1cm} $\underline{(2n-2-3i,1),(-(2n-4-3i),1),(2n-4-3i,1),(-(2n-6-3i),1),(2n-3-3i,1)},\ldots,(-\frac{3n+2}{2},1),$
\hspace*{1cm} $(\frac{3n+8}{2},1),(-\frac{3n+4}{2},1),(\frac{3n+4}{2},1),(-\frac{3n}{2},1),(\frac{3n+6}{2},1))$, $0 \leq i \leq \frac{n-12}{6}$;\\
$T_3=((\frac{3n}{2},1),(-\frac{3n-2}{2},1),(\frac{3n+2}{2},1))$;\\
$T_4=((-\frac{3n-6}{2},1),(\frac{3n-2}{2},1),(-\frac{3n-4}{2},1),(\frac{3n-6}{2},1),$
$(-\frac{3n-8}{2},1),(\frac{3n-4}{2},1),\ldots, \underline{(-(\frac{3n-6}{2}-3i),1),(\frac{3n-2}{2}-3i,1),}$\\
\hspace*{1cm} $ \underline{(-(\frac{3n-4}{2}-3i),1),(\frac{3n-6}{2}-3i,1),(-(\frac{3n-8}{2}-3i),1),(\frac{3n-4}{2}-3i,1)},\ldots,(-(n+3),1),(n+5,1),$\\
\hspace*{1cm}  $(-(n+4),1),(n+3,1),(-(n+2),1),(n+4,1))$, $0 \leq i \leq \frac{n-12}{6}$;\\
$T_5=((-n,1),(n+1,1),(n+1,3))$;\\
$T_6=((n+2,2),(-(n+4),2),(n+4,2),(-(n+3),2),(n+6,2),(-(n+5),2),\ldots,\underline{(n+2+3i,2),(-(n+4+3i),2),}$
 \hspace*{1cm} $\underline{(n+4+3i,2),(-(n+3+3i),2),(n+6+3i,2),(-(n+5+3i),2)},\ldots,(\frac{3n-8}{2},2),(-\frac{3n-4}{2},2),(\frac{3n-4}{2},2),$
 \hspace*{1cm} $(-\frac{3n-6}{2},2),(\frac{3n}{2},2),(-\frac{3n-2}{2},2))$, $0 \leq i \leq \frac{n-12}{6}$;\\
$T_7=((\frac{3n-2}{2},2),(-\frac{3n+2}{2},2),(\frac{3n+2}{2},2),(-\frac{3n}{2},2),(\frac{3n+8}{2},2),(-\frac{n-2}{2},2))$;\\
$T_8=((\frac{3n+4}{2},2),(-\frac{3n+8}{2},2),(\frac{3n+6}{2},2),(-\frac{3n+4}{2},2),(\frac{3n+14}{2},2),(-\frac{3n+6}{2},2),\ldots,\underline{(\frac{3n+4}{2}+3i,2),(-(\frac{3n+8}{2}+3i),2),}$
 \hspace*{1cm} $\underline{(\frac{3n+6}{2}+3i,2),(-(\frac{3n+4}{2}+3i),2),(\frac{3n+14}{2}+3i,2),(-(\frac{3n+6}{2}+3i),2)},\ldots,(2n-7,2),(-(2n-5),2),$\\
 \hspace*{1cm} $(2n-6,2),(-(2n-7),2),(2n-2,2),(-(2n-6),2))$, $0 \leq i \leq \frac{n-18}{6}$;\\
$T_9=((2n-4,2),(-(2n-3),2),(2n-3,2),(-(2n-2),2),(0,2),(2n-1,2),(-(2n-1),2),(2n,2),(-(2n-4),$
 \hspace*{1cm} $2),(2n+1,2),(-(n+2),2))$.
}

{\bf Case 1.2:} $k\equiv 11\pmod{24}$ and  $ k\geq 59$.\\
{\footnotesize
$S_6=((-\frac{n+2}{2},2),(\frac{n-8}{2},0),(-\frac{n}{2},2),(\frac{n-4}{2},0),(-\frac{n-8}{2},2),
(\frac{n-6}{2},0),\ldots,\underline{(-(\frac{n+2}{2}-3i),2),(\frac{n-8}{2}-3i,0),(-(\frac{n}{2}-3i),}$\\
\hspace*{1cm}  $\underline{2),(\frac{n-4}{2}-3i,0),(-(\frac{n-8}{2}-3i),2),(\frac{n-6}{2}-3i,0)},\ldots,(-8,2),(3,0),(-7,2),(5,0),(-3,2),(4,0))$,\\
\hspace*{1cm}  $0 \leq i \leq \frac{n-14}{6}$;\\
$S_7=((-5,2),(0,0),(-4,2),(2,0),(-1,2),(n-2,0),(-2,2))$.
}

The cycle $C_8$ is the concatenation of the sequences $\infty$,
$T_1$, $T_2$, $\ldots$, $T_9$, where\\
{\footnotesize
$T_1=((1,0),(-2n,2),(2n+1,1),(-2n,1),(2n,1),(-(2n-4),1),(-(2n+1),1),(-(2n-3),1),(2n-3,1),$
\hspace*{1cm} $(-(2n-2),1),(2n-1,1))$;\\
$T_2=((-(2n-7),1),(2n-2,1),(-(2n-5),1),(2n-6,1),(-(2n-6),1),(2n-4,1),\ldots,\underline{(-(2n-7-3i),1),}$\\
\hspace*{1cm} $ \underline{(2n-2-3i,1),(-(2n-5-3i),1),(2n-6-3i,1),(-(2n-6-3i),} \underline{1),(2n-4-3i,1)},\ldots,(-\frac{3n}{2},1),$\\
\hspace*{1cm} $(\frac{3n+10}{2},1),(-\frac{3n+4}{2},1),(\frac{3n+2}{2},1),(-\frac{3n+2}{2},1),(\frac{3n+6}{2},1))$, $0 \leq i \leq \frac{n-14}{6}$;\\
$T_3=((\frac{3n+4}{2},1),(-\frac{3n-4}{2},1))$;\\
$T_4=((\frac{3n}{2},1),(-\frac{3n-2}{2},1),(\frac{3n-4}{2},1),(-\frac{3n-6}{2},1),(\frac{3n-2}{2},1),$
$(-\frac{3n-10}{2},1),\ldots, \underline{(\frac{3n}{2}-3i,1),(-(\frac{3n-2}{2}-3i),1),}$\\
\hspace*{1cm} $\underline{(\frac{3n-4}{2}-3i,1),(-(\frac{3n-6}{2}-3i),1),(\frac{3n-2}{2}-3i,1),(-(\frac{3n-10}{2}-3i),1)},\ldots,(n+7,1),(-(n+6),1),$\\
\hspace*{1cm}  $(n+5,1),(-(n+4),1),(n+6,1),(-(n+2),1))$, $0 \leq i \leq \frac{n-14}{6}$;\\
$T_5=((n+3,1),(-(n+3),1),(n+4,1),(-n,1),(n+1,1),(n+1,3))$;\\
$T_6=((n+2,2),(-(n+4),2),(n+4,2),(-(n+3),2),(n+6,2),(-(n+5),2),\ldots,\underline{(n+2+3i,2),(-(n+4+3i),2),}$
 \hspace*{1cm} $\underline{(n+4+3i,2),(-(n+3+3i),2),(n+6+3i,2),(-(n+5+3i),2)},\ldots,(\frac{3n-4}{2},2),(-\frac{3n}{2},2),(\frac{3n}{2},2),$\\
 \hspace*{1cm} $(-\frac{3n-2}{2},2),(\frac{3n+4}{2},2),(-\frac{3n+2}{2},2))$, $0 \leq i \leq \frac{n-8}{6}$;\\
$T_7=((\frac{3n+2}{2},2),(-\frac{3n+6}{2},2),(\frac{3n+6}{2},2),(-\frac{3n+4}{2},2),(\frac{3n+10}{2},2),(-\frac{n-2}{2},2))$;\\
$T_8=((\frac{3n+12}{2},2),(-\frac{3n+10}{2},2),(\frac{3n+8}{2},2),(-\frac{3n+12}{2},2),(\frac{3n+16}{2},2),(-\frac{3n+8}{2},2),\ldots,\underline{(\frac{3n+12}{2}+3i,2),(-(\frac{3n+10}{2}+3i),}$
 \hspace*{1cm} $\underline{2),(\frac{3n+8}{2}+3i,2),(-(\frac{3n+12}{2}+3i),2),(\frac{3n+16}{2}+3i,2),(-(\frac{3n+8}{2}+3i),2)},\ldots,(2n-4,2),(-(2n-5),2),$
 \hspace*{1cm} $(2n-6,2),(-(2n-4),2),(2n-2,2),(-(2n-6),2))$, $0 \leq i \leq \frac{n-20}{6}$;\\
$T_9=((2n-1,2),(-(2n-2),2),(2n-3,2),(-(2n-1),2),(2n,2),(0,2),(2n+1,2),(-(2n-3),2),(-(n+2),2))$.
}

{\bf Case 1.3:} $k\equiv 19\pmod{24}$ and  $ k\geq 43$.\\
{\footnotesize
$S_6=((-\frac{n+2}{2},2),(\frac{n-8}{2},0),(-\frac{n}{2},2),(\frac{n-4}{2},0),(-\frac{n-8}{2},2),
(\frac{n-6}{2},0),\ldots,\underline{(-(\frac{n+2}{2}-3i),2),(\frac{n-8}{2}-3i,0),(-(\frac{n}{2}-3i),}$\\
\hspace*{1cm}  $\underline{2),(\frac{n-4}{2}-3i,0),(-(\frac{n-8}{2}-3i),2),(\frac{n-6}{2}-3i,0)},\ldots,(-9,2),(4,0),(-8,2),(6,0),(-4,2),(5,0))$,\\
\hspace*{1cm}  $0 \leq i \leq \frac{n-16}{6}$;\\
$S_7=((-6,2),(0,0),(-3,2),(2,0),(-5,2),(3,0),(-1,2),(n-2,0),(-2,2))$.
}

The cycle $C_8$ is the concatenation of the sequences $\infty$,
$T_1$, $T_2$, $\ldots$, $T_9$, where\\
{\footnotesize
$T_1=((1,0),(-2n,2),(2n+1,1),(2n-1,1),(-2n,1),(-(2n+1),1))$;\\
$T_2=((2n-4,1),(-(2n-3),1),(2n-2,1),(-(2n-4),1),(2n,1),(-(2n-2),1),\ldots,\underline{(2n-4-3i,1),}$\\
\hspace*{1cm} $ \underline{(-(2n-3-3i),1),(2n-2-3i,1),(-(2n-4-3i),1),} \underline{(2n-3i,1),(-(2n-2-3i),1)},\ldots,$\\
\hspace*{1cm} $(\frac{3n+2}{2},1),(-\frac{3n+4}{2},1),(\frac{3n+6}{2},1),(-\frac{3n+2}{2},1),(\frac{3n+10}{2},1),(-\frac{3n+6}{2},1))$, $0 \leq i \leq \frac{n-10}{6}$;\\
$T_3=((-\frac{3n}{2},1),(\frac{3n+4}{2},1))$;\\
$T_4=((-\frac{3n-4}{2},1),(\frac{3n}{2},1),(-\frac{3n-2}{2},1),(\frac{3n-4}{2},1),(-\frac{3n-6}{2},1),(\frac{3n-2}{2},1),\ldots, \underline{(-(\frac{3n-4}{2}-3i),1),(\frac{3n}{2}-3i,1),}$\\
\hspace*{1cm} $\underline{(-(\frac{3n-2}{2}-3i),1),(\frac{3n-4}{2}-3i,1),(-(\frac{3n-6}{2}-3i),1),(\frac{3n-2}{2}-3i,1)},\ldots,(-(n+3),1),(n+5,1),$\\
\hspace*{1cm}  $(-(n+4),1),(n+3,1),(-(n+2),1),(n+4,1))$, $0 \leq i \leq \frac{n-10}{6}$;\\
$T_5=((-n,1),(n+1,1),(n+1,3))$;\\
$T_6=((n+2,2),(-(n+4),2),(n+4,2),(-(n+3),2),(n+6,2),(-(n+5),2),\ldots,\underline{(n+2+3i,2),(-(n+4+3i),2),}$
 \hspace*{1cm} $\underline{(n+4+3i,2),(-(n+3+3i),2),(n+6+3i,2),(-(n+5+3i),2)},\ldots,(\frac{3n-6}{2},2),(-\frac{3n-2}{2},2),(\frac{3n-2}{2},2),$
 \hspace*{1cm} $(-\frac{3n-4}{2},2),(\frac{3n+2}{2},2),(-\frac{3n}{2},2))$, $0 \leq i \leq \frac{n-10}{6}$;\\
$T_7=((\frac{3n}{2},2),(-\frac{n-2}{2},2))$;\\
$T_8=((\frac{3n+8}{2},2),(-\frac{3n+6}{2},2),(\frac{3n+6}{2},2),(-\frac{3n+4}{2},2),(\frac{3n+4}{2},2),(-\frac{3n+2}{2},2),\ldots,\underline{(\frac{3n+8}{2}+3i,2),(-(\frac{3n+6}{2}+3i),2),}$
 \hspace*{1cm} $\underline{(\frac{3n+6}{2}+3i,2),(-(\frac{3n+4}{2}+3i),2),(\frac{3n+4}{2}+3i,2),(-(\frac{3n+2}{2}+3i),2)},\ldots,(2n-4,2),(-(2n-5),2),$\\
 \hspace*{1cm} $(2n-5,2),(-(2n-6),2),(2n-6,2),(-(2n-7),2))$, $0 \leq i \leq \frac{n-16}{6}$;\\
$T_9=((2n-1,2),(-(2n-2),2),(0,2),(2n+1,2),(-(2n-3),2),(2n-2,2),(-(2n-4),2),(2n-3,2),(-(2n-1),$
 \hspace*{1cm} $2),(2n,2),(-(n+2),2))$.
}

{\bf Case 2:}  $k\equiv 7\pmod{8}$ and $ k\geq 31$.

In this case the cycle $C_7$ is the concatenation of the sequences
$S_1$, $S_2$, $S_3$, $(n+3,2)$, $S_4$, $S_5$, $S_6$, and $S_7$, where\\
{\footnotesize
$S_1=((1,2),(-(2n-1),1))$;\\
$S_2=((2,2),(-2,0),\ldots,\underline{(2+i,2),(-(2+i),0)},\ldots,(\frac{n-1}{2},2),(-\frac{n-1}{2},0))$, $0 \leq i \leq \frac{n-5}{2}$;\\
$S_3=((\frac{n+3}{2},2),(-\frac{n+1}{2},0),\ldots,\underline{(\frac{n+3}{2}+i,2),(-(\frac{n+1}{2}+i),0)},\ldots,(n+1,2),(-n,0))$, $0 \leq i \leq \frac{n-1}{2}$;\\
$S_4=((n,0),(-(n-1),2),\ldots,\underline{(n-i,0),(-(n-1-i),2)},\ldots,(\frac{n+5}{2},0),(-\frac{n+3}{2},2))$, $0 \leq i \leq \frac{n-5}{2}$;\\
$S_5=((\frac{n+3}{2},0),(\frac{n+1}{2},2))$.
}

For the sequences $S_6$, $S_7$, and the cycle $C_8$, we distinguish 3 subcases.

{\bf Case 2.1:} $k\equiv 7\pmod{24}$ and   $ k\geq 31$. Here, $(a_2,b_2)=(\frac{3k-5}{4},2)$, $d_2=4$.\\
{\footnotesize
$S_6=((\frac{n-3}{2},0),(-\frac{n-5}{2},2),(\frac{n-1}{2},0),(-\frac{n-3}{2},2),(\frac{n+1}{2},0),(-\frac{n-1}{2},2),
\ldots,\underline{(\frac{n-3}{2}-3i,0),(-(\frac{n-5}{2}-3i),2),}$\\
\hspace*{1cm}  $\underline{(\frac{n-1}{2}-3i,0),(-(\frac{n-3}{2}-3i),2),(\frac{n+1}{2}-3i,0),(-(\frac{n-1}{2}-3i),2)},\ldots,(5,0),(-4,2),(6,0),(-5,2),$\\
\hspace*{1cm}  $(7,0),(-6,2))$, $0 \leq i \leq \frac{n-13}{6}$;\\
$S_7=((2,0),(0,2),(3,0),(-2,2),(4,0),(-3,2),(1,0),(-n,2),(-(n+1),0),(2n,2))$.
}

For $ k= 31, 55$, the cycle $C_8$ is listed as below.\\
{\footnotesize
$k=31$:\\
$C_8=(\infty, ( 0, 0),(-15, 2),(14, 1),(12, 1),(-14, 1),(-15, 1),(-12, 1),(15, 1),(-10, 1),(13, 1),(-11, 1),(10, 1),$ \\
\hspace*{0.9cm} $(-9, 1),(11, 1),(-7, 1),( 8, 1),( 8, 3),( 9, 2),(11, 2),(-11, 2),(15, 2),(-1, 2),(-13, 2),(-10, 2),(13, 2),$ \\
\hspace*{0.9cm} $(-4, 2),(-14, 2),(-8, 2),(12, 2),(-12, 2))$.
}\\
{\footnotesize
$k=55$:\\
$C_8=(\infty, ( 0, 0),(-27, 2),(26, 1),(-26, 1),(23, 1),(-27, 1),(27, 1),(-24, 1),(22, 1),(-23, 1),(25, 1),(-22, 1),$ \\
\hspace*{0.9cm} $(-20, 1),(24, 1),(-19, 1),(19, 1),(-18, 1),(21, 1),(-21, 1),(20, 1),(-16, 1),(16, 1),(-15, 1),(18, 1),$ \\
\hspace*{0.9cm} $(-17, 1),(17, 1),(-13, 1),(14, 1),(14, 3),(15, 2),(-18, 2),(17, 2),(-17, 2),(19, 2),(-19, 2),(18, 2),$ \\
\hspace*{0.9cm} $(-21, 2),(20, 2),(-20, 2),(22, 2),(-22, 2),(21, 2),(-25, 2),(27, 2),(-1, 2),(25, 2),(-7, 2),(24, 2),$ \\
\hspace*{0.9cm} $(-26, 2),(23, 2),(-24, 2),(-14, 2),(-16, 2),(-23, 2))$.
}

For any $ k\geq 79$, the cycle $C_8$ is the concatenation of the sequences $\infty$,
$T_1$, $T_2$, $\ldots$, $T_9$, where\\
{\footnotesize
$T_1=((0,0),(-(2n+1),2),(2n,1),(-2n,1),(2n-3,1)(-(2n+1),1),(2n+1,1)(-(2n-2),1))$;\\
$T_2=((2n-6,1),(-(2n-3),1),(2n-1,1),(-(2n-4),1),(2n-2,1),(-(2n-5),1),\ldots, \underline{(2n-6-3i,1),}$
\hspace*{1cm} $\underline{(-(2n-3-3i),1),} \underline{(2n-1-3i,1),(-(2n-4-3i),1),(2n-2-3i,1),(-(2n-5-3i),1)},\ldots,$\\
\hspace*{1cm} $(\frac{3n+7}{2},1),(-\frac{3n+13}{2},1),$ $(\frac{3n+17}{2},1),(-\frac{3n+11}{2},1),(\frac{3n+15}{2},1),(-\frac{3n+9}{2},1))$, $0\leq i \leq \frac{n-19}{6}$;\\
$T_3=((\frac{3n+5}{2},1),(-\frac{3n+7}{2},1),(\frac{3n+11}{2},1),(-\frac{3n+5}{2},1),(-\frac{3n+1}{2},1),(\frac{3n+9}{2},1),(-\frac{3n-1}{2},1),(\frac{3n-1}{2},1),(-\frac{3n-3}{2},1),$\\
\hspace*{1cm} $(\frac{3n+3}{2},1),(-\frac{3n+3}{2},1),(\frac{3n+1}{2},1),(-\frac{3n-7}{2},1))$;\\
$T_4=((\frac{3n-7}{2},1),(-\frac{3n-9}{2},1),(\frac{3n-3}{2},1),(-\frac{3n-5}{2},1),(\frac{3n-5}{2},1),(-\frac{3n-13}{2},1),\ldots,\underline{(\frac{3n-7}{2}-3i,1),(-(\frac{3n-9}{2}-3i),1),}$
\hspace*{1cm} $\underline{(\frac{3n-3}{2}-3i,1),(-(\frac{3n-5}{2}-3i),1),(\frac{3n-5}{2}-3i,1),(-(\frac{3n-13}{2}-3i),1)},\ldots,(n+3,1),(-(n+2),1),$\\
\hspace*{1cm} $(n+5,1),(-(n+4),1),(n+4,1),(-n,1))$, $0\leq i \leq \frac{n-13}{6}$;\\
$T_5=((n+1,1),(n+1,3))$;\\
$T_6=((n+2,2),(-(n+5),2),(n+4,2),(-(n+4),2),(n+6,2),(-(n+6),2),\ldots,\underline{(n+2+3i,2),(-(n+5+3i),2),}$
 \hspace*{1cm} $\underline{(n+4+3i,2),(-(n+4+3i),2),(n+6+3i,2),(-(n+6+3i),2)},\ldots,(\frac{3n-3}{2},2),(-\frac{3n+3}{2},2),(\frac{3n+1}{2},2),$
 \hspace*{1cm} $(-\frac{3n+1}{2},2),(\frac{3n+5}{2},2),(-\frac{3n+5}{2},2))$, $0\leq i \leq \frac{n-7}{6}$;\\
$T_7=((\frac{3n+3}{2},2),(\frac{3n+13}{2},2),(-\frac{3n+9}{2},2),(\frac{3n+7}{2},2),(-\frac{3n+11}{2},2),(\frac{3n+9}{2},2),(-\frac{n+1}{2},2))$;\\
$T_8=((\frac{3n+11}{2},2),(-\frac{3n+17}{2},2),(\frac{3n+15}{2},2),(-\frac{3n+15}{2},2),(\frac{3n+19}{2},2),(-\frac{3n+7}{2},2),\ldots,\underline{(\frac{3n+11}{2}+3i,2),(-(\frac{3n+17}{2}+3i),}$
 \hspace*{1cm} $\underline{2),(\frac{3n+15}{2}+3i,2),(-(\frac{3n+15}{2}+3i),2),(\frac{3n+19}{2}+3i,2),(-(\frac{3n+7}{2}+3i),2)},\ldots,(2n-7,2),(-(2n-4),2),$
 \hspace*{1cm}  $(2n-5,2),(-(2n-5),2),(2n-3,2),(-(2n-9),2))$, $0\leq i \leq \frac{n-25}{6}$;\\
$T_9=((2n-4,2),(-(2n-2),2),(2n-2,2),(-(2n-3),2),(-(n+1),2),(-(n+3),2),(-2n,2),(-(2n-6),2),$
 \hspace*{1cm} $(2n-1,2),(-1,2),(2n+1,2),(-(2n-1),2))$.
}

{\bf Case 2.2:} $k\equiv 15\pmod{24}$ and  $ k\geq 39$. Here, $(a_2,b_2)=(\frac{3k-5}{4},2)$, $d_2=4$. \\
{\footnotesize
$S_6=((\frac{n-3}{2},0),(-\frac{n-5}{2},2),(\frac{n-1}{2},0),(-\frac{n-3}{2},2),(\frac{n+1}{2},0),(-\frac{n-1}{2},2),
\ldots,\underline{(\frac{n-3}{2}-3i,0),(-(\frac{n-5}{2}-3i),2),}$\\
\hspace*{1cm}  $\underline{(\frac{n-1}{2}-3i,0),(-(\frac{n-3}{2}-3i),2),(\frac{n+1}{2}-3i,0),(-(\frac{n-1}{2}-3i),2)},\ldots,(6,0),(-5,2),(7,0),(-6,2),$\\
\hspace*{1cm}  $(8,0),(-7,2))$, $0 \leq i \leq \frac{n-15}{6}$;\\
$S_7=((3,0),(-3,2),(4,0),(-4,2),(5,0),(0,2),(2,0),(-2,2),(1,0),(-n,2),(-(n+1),0),(2n,2))$.
}

For $ k= 39$, , the cycle $C_8$ is listed as below.\\
{\footnotesize
$C_8=(\infty, ( 0, 0),(-19, 2),(18, 1),(-18, 1),(15, 1),(-19, 1),(19, 1),(-16, 1),(-14, 1),(16, 1),(-15, 1),(17, 1),$ \\
\hspace*{0.9cm} $(-12, 1),(12, 1),(-11, 1),(14, 1),(-13, 1),(13, 1),(-9, 1),(10, 1),(10, 3),(11, 2),(-14, 2),(13, 2), $ \\
\hspace*{0.9cm} $(-13, 2),(15, 2),(-15, 2),(14, 2),(-5, 2),(19, 2),(16, 2),(-16, 2),(17, 2),(-1, 2),(-17, 2),(-12, 2),$ \\
\hspace*{0.9cm} $(-10, 2),(-18, 2))$.
}

For any $ k\geq 63$, the cycle $C_8$ is the concatenation of the sequences $\infty$,
$T_1$, $T_2$, $\ldots$, $T_9$, where\\
{\footnotesize
$T_1=((0,0),(-(2n+1),2),(2n,1),(-2n,1),(2n-3,1)(-(2n+1),1),(2n+1,1)(-(2n-2),1))$;\\
$T_2=((2n-6,1),(-(2n-3),1),(2n-1,1),(-(2n-4),1),(2n-2,1),(-(2n-5),1),\ldots,$ $ \underline{(2n-6-3i,1),}$
\hspace*{1cm} $ \underline{(-(2n-3-3i),1),}\underline{(2n-1-3i,1),(-(2n-4-3i),1),(2n-2-3i,1),(-(2n-5-3i),1)},\ldots,$\\
\hspace*{1cm} $(\frac{3n+3}{2},1),(-\frac{3n+9}{2},1),(\frac{3n+13}{2},1),(-\frac{3n+7}{2},1),(\frac{3n+11}{2},1),(-\frac{3n+5}{2},1))$, $0 \leq i \leq \frac{n-15}{6}$;\\
$T_3=((-\frac{3n+1}{2},1),(\frac{3n+5}{2},1),(-\frac{3n+3}{2},1),(\frac{3n+7}{2},1),(-\frac{3n-3}{2},1))$;\\
$T_4=((\frac{3n-3}{2},1),(-\frac{3n-5}{2},1),(\frac{3n+1}{2},1),(-\frac{3n-1}{2},1),(\frac{3n-1}{2},1),(-\frac{3n-9}{2},1),\ldots,\underline{(\frac{3n-3}{2}-3i,1),(-(\frac{3n-5}{2}-3i),1),}$
\hspace*{1cm} $\underline{(\frac{3n+1}{2}-3i,1),(-(\frac{3n-1}{2}-3i),1),(\frac{3n-1}{2}-3i,1),(-(\frac{3n-9}{2}-3i),1)},\ldots,(n+3,1),(-(n+2),1),(n+5,1),$
\hspace*{1cm} $(-(n+4),1),(n+4,1),(-n,1))$, $0 \leq i \leq \frac{n-9}{6}$;\\
$T_5=((n+1,1),(n+1,3))$;\\
$T_6=((n+2,2),(-(n+5),2),(n+4,2),(-(n+4),2),(n+6,2),(-(n+6),2),\ldots,\underline{(n+2+3i,2),(-(n+5+3i),2),}$
 \hspace*{1cm} $\underline{(n+4+3i,2),(-(n+4+3i),2),(n+6+3i,2),(-(n+6+3i),2)},\ldots,(\frac{3n-5}{2},2),(-\frac{3n+1}{2},2),(\frac{3n-1}{2},2),$
 \hspace*{1cm} $(-\frac{3n-1}{2},2),(\frac{3n+3}{2},2),(-\frac{3n+3}{2},2))$, $0 \leq i \leq \frac{n-9}{6}$;\\
$T_7=((\frac{3n+1}{2},2),(-\frac{n+1}{2},2),(\frac{3n+11}{2},2),(\frac{3n+5}{2},2))$;\\
$T_8=((-\frac{3n+7}{2},2),(\frac{3n+17}{2},2),(-\frac{3n+5}{2},2),(\frac{3n+9}{2},2),(-\frac{3n+9}{2},2),(\frac{3n+7}{2},2),\ldots,\underline{(-(\frac{3n+7}{2}+3i),2),(\frac{3n+17}{2}+3i,2),}$
 \hspace*{1cm} $\underline{(-(\frac{3n+5}{2}+3i),2),(\frac{3n+9}{2}+3i,2),(-(\frac{3n+9}{2}+3i),2),(\frac{3n+7}{2}+3i,2)},\ldots,(-(2n-7),2),(2n-2,2),$\\
 \hspace*{1cm}  $(-(2n-8),2),(2n-6,2),(-(2n-6),2),(2n-7,2))$, $0 \leq i \leq \frac{n-21}{6}$;\\
$T_9=((-(2n-4),2),(2n-4,2),(-(2n-3),2),(2n-3,2),(-(2n-2),2),(2n-1,2),(-1,2),(-(2n-1),2),$
 \hspace*{1cm} $(-(n+1),2),(-(n+3),2),(2n+1,2),(-(2n-5),2),(-2n,2))$.
}

{\bf Case 2.3:} $k\equiv 23\pmod{24}$ and  $ k\geq 47$. Here, $(a_2,b_2)=(\frac{3k-1}{4},2)$, $d_2=2$. \\
{\footnotesize
$S_6=((\frac{n-3}{2},0),(-\frac{n-5}{2},2),(\frac{n-1}{2},0),(-\frac{n-3}{2},2),(\frac{n+1}{2},0),(-\frac{n-1}{2},2),
\ldots,\underline{(\frac{n-3}{2}-3i,0),(-(\frac{n-5}{2}-3i),2),}$\\
\hspace*{1cm}  $\underline{(\frac{n-1}{2}-3i,0),(-(\frac{n-3}{2}-3i),2),(\frac{n+1}{2}-3i,0),(-(\frac{n-1}{2}-3i),2)},\ldots,(7,0),(-6,2),(8,0),(-7,2),(9,0),$\\
\hspace*{1cm}  $(-8,2))$, $0 \leq i \leq \frac{n-17}{6}$;\\
$S_7=((4,0),(-3,2),(3,0),(-5,2),(6,0),(-4,2),(5,0),(0,2),(2,0),(-2,2),(1,0),(-n,2),(-(n+1),0),(2n,2))$.
}

For $ k= 47$, the cycle $C_8$ is listed as below. \\
{\footnotesize
$C_8=(\infty, ( 0, 0),(-23, 2),(22, 1),(-22, 1),(19, 1),(-23, 1),(23, 1),(-20, 1),(20, 1),(-19, 1),(-17, 1),(21, 1),$ \\
\hspace*{0.9cm} $(-16, 1),(17, 1),(-18, 1),(18, 1),(-14, 1),(14, 1),(-13, 1),(16, 1),(-15, 1),(15, 1),(-11, 1),(12, 1),$ \\
\hspace*{0.9cm} $(12, 3),(13, 2),(-16, 2),(18, 2),(-15, 2),(17, 2),(-14, 2),(16, 2),(19, 2),(-6, 2),(21, 2),(-21, 2),(15, 2),$ \\
\hspace*{0.9cm} $(-20, 2),(-1, 2),(23, 2),(-17, 2),(20, 2),(-19, 2),(-13, 2),(-22, 2),(-18, 2))$.
}

For any $ k\geq 71$, the cycle $C_8$ is the concatenation of the sequences $\infty$,
$T_1$, $T_2$, $\ldots$, $T_9$, where\\
{\footnotesize
$T_1=((0,0),(-(2n+1),2),(2n,1),(-2n,1),(2n-3,1)(-(2n+1),1),(2n+1,1)(-(2n-2),1))$;\\
$T_2=((2n-6,1),(-(2n-3),1),(2n-1,1),(-(2n-4),1),(2n-2,1),(-(2n-5),1),\ldots,$ $ \underline{(2n-6-3i,1),}$
\hspace*{1cm} $ \underline{(-(2n-3-3i),1),}\underline{(2n-1-3i,1),(-(2n-4-3i),1),(2n-2-3i,1),(-(2n-5-3i),1)},\ldots,$\\
\hspace*{1cm} $(\frac{3n+5}{2},1),(-\frac{3n+11}{2},1),(\frac{3n+15}{2},1),(-\frac{3n+9}{2},1),(\frac{3n+13}{2},1),(-\frac{3n+7}{2},1))$, $0 \leq i \leq \frac{n-17}{6}$;\\
$T_3=((\frac{3n+7}{2},1),(-\frac{3n+5}{2},1),(-\frac{3n+1}{2},1),(\frac{3n+9}{2},1),(-\frac{3n-1}{2},1),(\frac{3n+1}{2},1),(-\frac{3n+3}{2},1),(\frac{3n+3}{2},1),(-\frac{3n-5}{2},1))$;\\
$T_4=((\frac{3n-5}{2},1),(-\frac{3n-7}{2},1),(\frac{3n-1}{2},1),(-\frac{3n-3}{2},1),(\frac{3n-3}{2},1),(-\frac{3n-11}{2},1),\ldots,\underline{(\frac{3n-5}{2}-3i,1),(-(\frac{3n-7}{2}-3i),1),}$
\hspace*{1cm} $\underline{(\frac{3n-1}{2}-3i,1),(-(\frac{3n-3}{2}-3i),1),(\frac{3n-3}{2}-3i,1),(-(\frac{3n-11}{2}-3i),1)},\ldots,(n+3,1),(-(n+2),1),$\\
\hspace*{1cm} $(n+5,1),(-(n+4),1),(n+4,1),(-n,1))$, $0 \leq i \leq \frac{n-11}{6}$;\\
$T_5=((n+1,1),(n+1,3))$;\\
$T_6=((n+2,2),(-(n+5),2),(n+7,2),(-(n+4),2),(n+6,2),(-(n+3),2),\ldots,\underline{(n+2+3i,2),(-(n+5+3i),2),}$
 \hspace*{1cm} $\underline{(n+7+3i,2),(-(n+4+3i),2),(n+6+3i,2),(-(n+3+3i),2)},\ldots,(\frac{3n-7}{2},2),(-\frac{3n-1}{2},2),(\frac{3n+3}{2},2),$
 \hspace*{1cm} $(-\frac{3n-3}{2},2),(\frac{3n+1}{2},2),(-\frac{3n-5}{2},2))$, $0 \leq i \leq \frac{n-11}{6}$;\\
$T_7=((\frac{3n-1}{2},2),(\frac{3n+5}{2},2),(-\frac{n+1}{2},2),(-\frac{3n+1}{2},2),(\frac{3n+7}{2},2),(-\frac{3n+5}{2},2),(\frac{3n+11}{2},2),(-\frac{3n+3}{2},2))$;\\
$T_8=((-\frac{3n+11}{2},2),(\frac{3n+17}{2},2),(-\frac{3n+7}{2},2),(\frac{3n+13}{2},2),(-\frac{3n+9}{2},2),(\frac{3n+9}{2},2),\ldots,\underline{(-(\frac{3n+11}{2}+3i),2),(\frac{3n+17}{2}+3i,2),}$
 \hspace*{1cm} $\underline{(-(\frac{3n+7}{2}+3i),2),(\frac{3n+13}{2}+3i,2),(-(\frac{3n+9}{2}+3i),2),(\frac{3n+9}{2}+3i,2)},\ldots,(-(2n-6),2),(2n-3,2),$\\
 \hspace*{1cm}  $(-(2n-8),2),(2n-5,2),(-(2n-7),2),(2n-7,2))$, $0 \leq i \leq \frac{n-23}{6}$;\\
$T_9=((-(2n-3),2),(2n-2,2),(-(2n-5),2),(2n-1,2),(-(2n-1),2),(-1,2),(2n+1,2),(-(2n-4),2),$
 \hspace*{1cm} $(2n-4,2),(-2n,2),(-(n+2),2),(n+4,2),(-(2n-2),2))$.
}
\qed

\section{Concluding remarks}

Combining Theorems~\ref{3-14}, \ref{2kt+1}, and Lemmas~\ref{4k+1}, \ref{4k+3}, we have proved Theorem~\ref{main}. This completes the proof of the existence of
almost resolvable cycle systems with odd cycle length. For the even case,  Lemma~\ref{A} is still useful for certain subcases.

We are working on the case $k=4m+2$.  But it should be mentioned that Lemma~\ref{A} can not be applied to solve the  existence of a $k$-almost resolvable cycle system when $k=2^m$. So there is still a long way to go before the whole problem can be solved completely.

As an application, almost resolvable cycle systems can be used to construct some
solutions to the Hamilton-Waterloo problem \cite{WLC}. The Hamilton-Waterloo problem HWP$(v;m,n;\alpha,\beta)$ is the problem of
determining whether the complete graph $K_v$ (for $v$ odd) or $K_v$ minus a
$1$-factor $I$ (for $v$ even) has a $2$-factorization in which there
are exactly $\alpha$ $C_m$-factors and $\beta$ $C_n$-factors.
We denote by HWP$(v;m,n)$ the set of $(\alpha,\beta)$ for which a solution to
HWP$(v;m,n;\alpha,\beta)$ exists.
For recent results on the Hamilton-Waterloo problem,
we refer the reader to \cite{AKK,BB,BD2,BDT,BDT2,BDT3,KP,LF,MT,OO,WC,WCC}. As a by-product, the following theorem can be obtained by combining Theorem~\ref{main} of this paper, and Theorems 1.4 and 3.5, Constructions 3.11 and 3.13 in \cite{WLC}.

\begin{theorem}
\label{Q2}

If $k\geq 3$ is odd and $t\geq 1$, then $(\alpha,\beta)\in$ {\rm HWP}$ (k(2kt+1); k, 2kt+1)$ if and only if $\alpha, \beta \geq 0$ and $\alpha+\beta=\frac{k(2kt+1)-1}{2}$,
 except possibly when:

$1.$ $t=1.$

$\ \ \ \textcircled{\small{1}}$ $k=5:$ $\beta \in \{1, 2, 3\}$;

$\ \ \ \textcircled{\small{2}}$ $k=7:$ $\beta \in \{ 1, 2, 3, 5\}$;

$\ \ \ \textcircled{\small{3}}$ $k\geq 9:$ $\beta \in \{ 1, 2, 3, 5, 7\}$.

$2.$ $t=2.$

$\ \ \ \textcircled{\small{1}}$ $k=3:$ $\beta \in [1, 2k-1]\cup \{2k+1, 2k+3\}$;

$\ \ \ \textcircled{\small{2}}$ $k=5, 7, 9 :$ $\beta \in \{ 1, 2, 3, 5, 7\}$.

$3.$ $t\geq 3.$

$\ \ \ \textcircled{\small{1}}$ $k=3:$

$\ \ \ \  \bullet$ $t$ is odd $:$ $\beta \in \{1,3,5,\dots,3t-4, 3t-2\} \cup \{2, 9t-3, 9t-1\}$

$\ \ \ \  \bullet$ $t$ is even $:$ $\beta \in \{1,3,5,\dots,3t+1, 3t+3\} \cup \{2, 9t-3, 9t-1\}$;

$\ \ \ \textcircled{\small{2}}$ $k\geq 5:$ $\beta \in \{ 1, 2, 3, 5, 7\}$.
\end{theorem}

\vspace{5pt}

\noindent {\bf Acknowledgments}

We would like to thank the anonymous referees for their careful read
and helpful comments and suggestions which greatly improved the quality of this paper.

\end{document}